\documentclass[12pt]{amsart}

\usepackage{graphicx}
\usepackage{amsmath}
\usepackage{hyperref}

\usepackage{tikz}

\def\frak{\mathfrak}

\def\R{\mathbb{R}}

\def\cD{\mathcal{D}}

\def\al{\alpha}

\def\ga{\gamma}
\def\de{\delta}

\def\om{\omega}
\def\Ga{\Gamma}
\def\De{\Delta}

\def\Om{\Omega}

\newcommand{\der}{{\rm d}}

\numberwithin{equation}{section}

\newtheorem{theorem}{Theorem}[section]

\newtheorem{proposition}[theorem]{Proposition}
\newtheorem{corollary}[theorem]{Corollary}
\theoremstyle{remark}

\theoremstyle{remark}

\usepackage{amssymb,stmaryrd}
\usepackage{amscd}

\newcommand{\qr}{
\begin{center}
\includegraphics[scale=0.5]{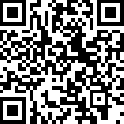}\\

\text{Scan the QR code to view more articles from the author}
\end{center}}

\makeatletter
\@namedef{subjclassname@2020}{%
  \textup{2020} Mathematics Subject Classification}
\makeatother

\author{Matthew Randall}
\address{Institute of Mathematical Sciences \\
ShanghaiTech University\\
393 Middle Huaxia Road\\
Shanghai 201210\\
China}
\email{mjrandall@shanghaitech.edu.cn}

\title{Local equivalence of some maximally symmetric $(2,3,5)$-distributions II}


\begin{document}

\begin{abstract}
We show the change of coordinates that maps the maximally symmetric $(2,3,5)$-distribution given by solutions to the $k=\frac{2}{3}$ and $k=\frac{3}{2}$ generalised Chazy equation to the flat Cartan distribution. This establishes the local equivalence between the maximally symmetric $k=\frac{2}{3}$ and $k=\frac{3}{2}$ generalised Chazy distribution and the flat Cartan or Hilbert-Cartan distribution. We give the set of vector fields parametrised by solutions to the $k=\frac{2}{3}$ and $k=\frac{3}{2}$ generalised Chazy equation and the corresponding Ricci-flat conformal scale that bracket-generate to give the split real form of $\frak{g}_2$.
\end{abstract}

\maketitle

\pagestyle{myheadings}
\markboth{Randall}{Local equivalence of some maximally symmetric $(2,3,5)$-distributions II}

\tableofcontents

\section{Introduction}

Let $\cD$ be a maximally non-integrable rank 2 distribution on a 5-manifold $M$. The maximally non-integrable condition of $\cD$ determines a filtration of the tangent bundle $TM$ given by
\[
\cD \subset [\cD,\cD] \subset [\cD,[\cD,\cD]]\cong TM.
\]
The distribution $[\cD, \cD]$ has rank 3 while the full tangent space $TM$ has rank 5, hence such a geometry is also known as a $(2,3,5)$-distribution. Let $M_{xyzpq}$ denote the 5-dimensional mixed order jet space $J^{2,0}(\R,\R^2) \cong J^2(\R,\R)\times \R$ with local coordinates given by $(x,y,z,p,q)=(x,y,z,y',y'')$ (see also \cite{tw13}, \cite{tw14}). Let $\cD_{\varphi(x,y,z,y',y'')}$ denote the maximally non-integrable rank 2 distribution on $M_{xyzpq}$ associated to the underdetermined differential equation $z'=\varphi(x,y,z,y',y'')$. This means that the distribution is annihilated by the following three 1-forms
\begin{align*}
\om_1=\der y-p \der x, \qquad \om_2=\der p-q \der x, \qquad  \om_3=\der z-\varphi(x,y,z,p,q) \der x.
\end{align*}
Such a distribution $\cD_{\varphi(x,y,z,y',y'')}$ is said to be in Monge normal form (see page 90 of \cite{tw13}). 
The historically important example is the 1-forms associated to the Hilbert-Cartan distribution obtained when $\varphi(x,y,z,p,q)=q^2$ \cite{cartan1910}. This distribution gives the flat model of a $(2,3,5)$-distribution and is associated to the Hilbert-Cartan equation $z'=(y'')^2$ (see Section 5 of \cite{conf} for a discussion of this equation).

In Section 5 of \cite{conf}, it is shown how to associate canonically to such a $(2,3,5)$-distribution a conformal class of metrics of split signature $(2,3)$ (henceforth known as Nurowski's conformal structure or Nurowski's conformal metrics) such that the rank 2 distribution is isotropic with respect to any metric in the conformal class. The method of equivalence \cite{cartan1910} (also see the introduction to \cite{annur}, Section 5 of \cite{conf} and \cite{Strazzullo}) produces the 1-forms $(\theta_1, \theta_2,\theta_3, \theta_4, \theta_5)$ that give a coframing for Nurowski's metric. These 1-forms satisfy the structure equations
\begin{align*}
\der \theta_1&=\theta_1\wedge (2\Om_1+\Om_4)+\theta_2\wedge \Om_2+\theta_3 \wedge \theta_4,\nonumber\\
\der \theta_2&=\theta_1\wedge\Om_3+\theta_2\wedge (\Om_1+2\Om_4)+\theta_3 \wedge \theta_5,\nonumber\\
\der \theta_3&=\theta_1\wedge\Om_5+\theta_2\wedge\Om_6+\theta_3\wedge (\Om_1+\Om_4)+\theta_4 \wedge \theta_5,\\
\der \theta_4&=\theta_1\wedge\Om_7+\frac{4}{3}\theta_3\wedge\Om_6+\theta_4\wedge \Om_1+\theta_5 \wedge \Om_2,\nonumber\\
\der \theta_5&=\theta_2\wedge \Om_7-\frac{4}{3}\theta_3\wedge \Om_5+\theta_4\wedge\Om_3+\theta_5\wedge \Om_4,\nonumber
\end{align*}
where $(\Om_1, \ldots, \Om_7)$ and two additional 1-forms $(\Om_8, \Om_9)$ together define a rank 14 principal bundle over the 5-manifold $M$ (see \cite{cartan1910} and Section 5 of \cite{conf}). A representative metric in Nurowski's conformal class \cite{conf} is given by
\begin{align}\label{metric}
g=2 \theta_1 \theta_5-2\theta_2 \theta_4+\frac{4}{3}\theta_3 \theta_3.
\end{align}
When $g$ has vanishing Weyl tensor, the distribution is called maximally symmetric and has split $G_2$ as its group of local symmetries. For further details about the curvature invariant, see the introduction to \cite{annur} and Section 5 of \cite{conf}. For further discussion on the relationship between maximally symmetric $(2,3,5)$-distributions and the automorphism group of the split octonions, see Section 2 of \cite{tw13}.

For example, when $\varphi(x,y,z,p,q)=q^m$, we obtain the distribution associated to the equation $z'=(y'')^m$. For such distributions, Nurowski's metric \cite{conf} given by (\ref{metric}) has vanishing Weyl tensor precisely when $m \in \{-1,\frac{1}{3},\frac{2}{3},2\}$. For the values of $m=-1,\frac{1}{3}$ and $\frac{2}{3}$ these maximally symmetric distributions are all locally diffeomorphic to the $m=2$ Hilbert-Cartan case. This means that for the distributions of the form $\varphi(x,y,z,p,q)=q^m$, with $m \in \{-1,\frac{1}{3},\frac{2}{3}\}$, we can redefine local coordinates to obtain the Hilbert-Cartan distribution (see for instance \cite{DK}).

$(2,3,5)$-distributions also arise from the study of the configuration space of two surfaces rolling without slipping or twisting over each other \cite{AN14}, \cite{BH} and \cite{BM}. The configuration space can be realised as the An-Nurowski circle twistor distribution \cite{AN14} and in the case of two spheres with radii in the ratio $1:3$ rolling without slipping or twisting over each other, there is again maximal $G_2$ symmetry. 

In the work of \cite{r17}, a description of maximally symmetric $(2,3,5)$-distributions obtained from Pfaffian systems with $SU(2)$ symmetry was discussed and its relationship with the rolling distribution was investigated. In particular, the An-Nurowski circle twistor bundle can be realised by considering the Riemannian surface element of the unit sphere arising from one copy of $SU(2)$ and the other Riemannian surface element with Gaussian curvature $9$ or $\frac{1}{9}$ from another copy of $SU(2)$. Both Lie algebras of $su(2)$ are parametrised by the left-invariant vector fields. See \cite{r17} for further details. 

In \cite{r21a}, using a parametrisation of the Lie algebra of $sl_2$ given by the second prolongation of the group of unimodular fractional linear transformations, a Monge normal form for the rolling distribution in the case of hyperboloid surfaces is obtained. Using this, the author in \cite{r21b} was able to derive the analogous parametrisation for a complexified $su(2)$ in the sphere rolling distribution and find a change of coordinates that bring it into the Monge normal form with 
\[
\varphi(x,y,z,p,q)=qz^2-\frac{1}{\al^2+1}(\sqrt{q}z-\frac{1}{2\sqrt{q}x})^2. 
\]
Here $\al$ is a complex number, and the maximally symmetric case is obtained 
whenever $\al^2=-\frac{1}{9}$ or $\al^2=-9$. 

In the aforementioned maximally symmetric case, we found in \cite{r21b} the change of coordinates that maps the rolling distribution into the flat Cartan distribution, and therefore as a corollary into the Hilbert-Cartan distribution. This establishes the local equivalence between the maximally symmetric rolling model and the flat Cartan or Hilbert-Cartan distribution. The vector fields that bracket-generate to give the split real form of the Lie algebra of $\frak{g}_2$ were also found, with two of the vector fields in the bracket-generating set given by the span of the rolling distribution. They were presented in Theorems 4.1 and 4.2 of \cite{r21b}.

In this article we continue our investigations into the maximally symmetric $(2,3,5)$-distributions given by solutions to the $k=\frac{2}{3}$ and $k=\frac{3}{2}$ generalised Chazy equation. We consider distributions of the form $\varphi(x,y,z,p,q)=\frac{q^2}{H''(x)}$. The Weyl tensor vanishes in the case where $H(x)$ satisfies the 6th-order ordinary differential equation (ODE) known as Noth's equation \cite{annur}. The 6th-order ODE can be solved by the generalised Chazy equation with parameter $k=\frac{3}{2}$ and its Legendre dual is another 6th-order ODE that can be solved by the generalised Chazy equation with parameter $k=\frac{2}{3}$ \cite{r16}.  

For such maximally symmetric distributions described locally by a certain function $\varphi(x,q)=\frac{q^2}{H''(x)}$ where $H(x)$ satisfies Noth's equation, it was found in \cite{r19} the corresponding Ricci-flat representatives in Nurowski's conformal class. This involves solving a second-order differential equation (see equations (1.2) and (1.5) of \cite{r19}) to find the conformal scale in which the Ricci tensor of the conformally rescaled metric vanishes, which turns out to be related to the solutions of Noth's equation. The second-order differential equation that determines the conformal scale for Ricci-flatness involves solutions of the generalised Chazy equation with parameter $k=3$ and in the dual case $k=2$. These are the results of Theorems 3.1 and 3.2 of \cite{r19}.

Once the Ricci-flat representatives are found, the metric is both Ricci-flat and conformally flat and the challenge is to redefine local coordinates to obtain the flat Cartan distribution. This can be quite easy to find for some distributions such as those of the form $\varphi(x,y,z,p,q)=q^m$, with $m \in \{-1,\frac{1}{3},\frac{2}{3}\}$, but can also take a longer time to obtain for the ones that we consider here. 

In this article we find the change of coordinates that maps the maximally symmetric generalised Chazy distribution into the flat Cartan distribution, and therefore as a corollary into the Hilbert-Cartan distribution. This establishes the local equivalence between the maximally symmetric generalised Chazy model and the flat Cartan or Hilbert-Cartan distribution. For the maximally symmetric generalised Chazy distribution, we write down the vector fields that bracket-generate to give the split real form of the Lie algebra of $\frak{g}_2$, with two of the vector fields in the bracket-generating set given by the span of the maximally symmetric generalised Chazy distribution. These are presented in Theorems \ref{g2a} and \ref{g2b}. They depend on solutions of the generalised Chazy equation with parameters $k=\frac{2}{3}$ and $k=\frac{3}{2}$ and their corresponding Ricci-flat conformal scale. We then give as corollaries to the theorems in Corollary \ref{cora} and \ref{corb}, the examples when we take the solutions of the spin $\frac{3}{2}$ Lam\'e equation and solutions of the spin $4$ Lam\'e equation. 

This article can be viewed as a sequel to \cite{r21b}, both common in the purpose of writing down the vector fields associated to maximally symmetric $(2,3,5)$-distributions that bracket generate split $\frak{g}_2$. It can also be seen as a sequel to \cite{r19}, since the technical details are similar and we make use of the results about the conformal factor for Ricci-flatness derived there. Nonetheless, we also try to make this article as self-contained as possible, since the theme of integrable differential equations appear here and also so that readers do not have to rely too much on cross references. 

The computations here are done utilising heavily the \texttt{DifferentialGeometry} package in MAPLE 2018. 

\section{Flat Cartan distribution}
This section is reproduced from Section 2 in \cite{r21b}, in order to make the article self-contained. We recall that the coframe data of the canonically maximally symmetric $(2,3,5)$-distribution is given by 
\begin{align*}
\der\theta_1=\theta_3 \wedge \theta_4, \quad
\der\theta_2=\theta_3 \wedge \theta_5,\quad
\der\theta_3=\theta_4 \wedge \theta_5,\quad
\der\theta_4=0,\quad
\der\theta_5=0.
\end{align*}
This is the historic case studied by Cartan (\cite{cartan1893}, \cite{cartan1910}) and Engel (\cite{Engel}, \cite{Engel2}). 
There are local coordinates $(a_1,a_2,a_3,a_4,a_5)$ (see pages 159--160 of \cite{cartan1910}) such that 
\begin{align}\label{t1}
\theta_1&=\der a_1+\left(a_3+\frac{1}{2}a_4 a_5\right)\der a_4,\\
\label{t2}
\theta_2&=\der a_2+\left(a_3-\frac{1}{2}a_4 a_5\right)\der a_5,\\
\label{t3}
\theta_3&=\der a_3+\frac{1}{2}a_4 \der a_5-\frac{1}{2}a_5 \der a_4,\\
\theta_4&=\der a_4,\nonumber\\
\theta_5&=\der a_5.\nonumber
\end{align}
For non-zero constant $k$, the weighted rescaling
\[
(\al_1,\al_2,\al_3, \al_4, \al_5)\mapsto (k^3\al_1,k^3\al_2,k^2 \al_3, k\al_4, k\al_5)
\]
preserves $\theta_1$, $\theta_2$ and $\theta_3$.
We shall refer to this distribution annhilated by the 1-forms $\{\theta_1, \theta_2, \theta_3\}$ in (\ref{t1}), (\ref{t2}), (\ref{t3}) as the flat Cartan distribution. To facilitate our writing of the vector fields that bracket-generate the Lie algebra of split ${\frak g}_2$, let us pass to the 1-forms 
\begin{align*}
\Theta_1&=dc_1-2c_4dc_3-4c_3dc_4,\\
\Theta_2&=dc_2+2c_5dc_3+4c_3 dc_5,\\
\Theta_3&=dc_3+c_5dc_4-c_4dc_5,
\end{align*}
by taking the change of coordinates
\begin{align*}
(c_1,c_2,c_3,c_4,c_5)=\left(6a_1-2a_3a_4+a_4^2a_5,6a_2-2a_3a_5-a_4a_5^2,2a_3,-a_4,a_5\right).
\end{align*}
It follows that
\begin{align*}
\Theta_1=6\theta_1+2a_4 \theta_3,\quad \Theta_2=6\theta_2+2a_5\theta_3,\quad \Theta_3=2\theta_3,
\end{align*}
so the 1-forms $\{\Theta_1, \Theta_2, \Theta_3\}$ are in the span of $\{\theta_1,\theta_2,\theta_3\}$.
If we take $\frak{r_1}=c_5$,$\frak{r_2}=c_4$, $\frak{r_3}=c_3$, $\frak{r_4}=\frac{1}{2}(c_2+3c_3c_5)$, $\frak{r_5}=\frac{1}{2}(c_1-3c_3c_4)$, then the 1-forms 
\begin{align*}
\der\frak{r_3}+\frak{r_1}\der \frak{r_2}-\frak{r_2}\der \frak{r_1}&=\Theta_3,\\
\der\frak{r_4}+\frac{1}{2}(\frak{r_3}\der \frak{r_1}-\frak{r_1}\der \frak{r_3})&=\frac{1}{2}\Theta_2,\\
\der\frak{r_5}+\frac{1}{2}(\frak{r_2}\der \frak{r_3}-\frak{r_3}\der \frak{r_2})&=\frac{1}{2}\Theta_1,
\end{align*}
obtained are those annihilating the flat Engel distribution as given in \cite{Engel} and \cite{Engel2}.

Let us write down the vector fields
\begin{align*}
Z^1=\partial_{c_3}+2c_5\partial_{c_2}-2c_4\partial_{c_1},\\
Z^2=\partial_{c_4}+4c_3\partial_{c_1}-2c_5\partial_{c_3},\\
Z^3=\partial_{c_5}+2c_4\partial_{c_3}-4c_3\partial_{c_2},
\end{align*}
and define
\begin{align}\label{cdist}
S^1=Z^2+c_5Z^1,\hspace{12pt} S^2=Z^3-c_4Z^1, \hspace{12pt} S^3=-c_1Z^2+c_2Z^3-(c_1c_5+c_2c_4+c_3^2)Z^1.
\end{align}
The vector fields $S^1$ and $S^2$ are in the span of the distribution and are annihilated by the 1-forms $\{\Theta_1, \Theta_2, \Theta_3\}$. 
We say that the vector fields $\{S^1,S^2,S^3\}$ pairwise bracket-generate the Lie algebra of split ${\frak g}_2$ if the following holds:
defining
\begin{align*}
S^4&=[S^1,S^2],\quad S^5=[S^2,S^3],\quad S^6=[S^3,S^1],\\
L^1&=[S^1,S^4],\quad L^3=[S^2,S^5],\quad L^5=[S^3,S^6],\\
L^2&=[S^2,S^4],\quad L^4=[S^3,S^5],\quad L^6=[S^1,S^6]
\end{align*}
and
\begin{align*}
H=[S^2,S^6], \quad h=[S^4,S^3],
\end{align*}
we require that the set of vector fields
\begin{align*}
\{S^1,S^2,S^3,S^4,S^5,S^6,\frac{1}{4}(h-H),\frac{\sqrt{3}}{12}(h+H),L^1,L^2,L^3,L^4,L^5,L^6\}
\end{align*}
form the 14-dimensional Lie algebra of split ${\frak g}_2$ with the Cartan subalgebra spanned by $\frac{1}{4}(h-H)$ and $\frac{\sqrt{3}}{12}(h+H)$ and the root diagram given by the picture below with respect to this choice of the Cartan subalgebra. Further details about the formulas for the remaining vector fields and commutator relations can be found in the Appendix. 
\begin{figure}[h!]
\begin{tikzpicture}
	\draw [stealth-stealth](-1,0) -- (1,0);
\draw (1,0) node[anchor=west] {{\tiny $S^1$}};
\draw (-1,0) node[anchor=east] {{\tiny $S^5$}};
	\draw [stealth-stealth](0,-1.732) -- (0,1.732);
\draw (0,-1.732) node[anchor=north] {{\tiny $L^5$}};
\draw (0,1.732) node[anchor=south] {{\tiny $L^2$}};
\draw [stealth-stealth](-0.5,-0.866) -- (0.5,0.866);
\draw (-0.5,-0.866) node[anchor=north] {{\tiny $S^{3}$}};
\draw (0.5,0.866) node[anchor=south] {{\tiny $S^4$}};
\draw [stealth-stealth](-1.5,-0.866) -- (1.5,0.866);
\draw (-1.5,-0.866) node[anchor=north] {{\tiny $L^4$}};
\draw (1.5,0.866) node[anchor=south] {{\tiny $L^1$}};
\draw [stealth-stealth](1.5,-0.866) -- (-1.5,0.866);
\draw(1.5,-0.866) node[anchor=north] {{\tiny $L^6$}};
\draw (-1.5,0.866) node[anchor=south] {{\tiny $L^3$}};
\draw [stealth-stealth](0.5,-0.866) -- (-0.5,0.866);
\draw (0.5,-0.866) node[anchor=north] {{\tiny $S^{6}$}};
\draw (-0.5,0.866) node[anchor=south] {{\tiny $S^2$}};
\end{tikzpicture}
\end{figure}

\begin{proposition}
The vector fields given in (\ref{cdist}) pairwise bracket-generate the Lie algebra of split $\frak{g}_2$.
\end{proposition}
Since $S^1$ and $S^2$ are spanned by the distribution, which is already given as part of the data, the non-trivial part in determining the generating set of the Lie algebra of split $\frak{g}_2$ for a maximally symmetric $(2,3,5)$-distribution is to find $S^3$. The way to find it is outlined as above. We find the change of coordinates that bring the 1-forms annihilating the distribution to the span of $\theta_1$, $\theta_2$, $\theta_3$ in the flat Cartan distribution. Then we determine the functions $(c_1,c_2,c_3,c_4,c_5)$ and write down the vector fields $Z^1$, $Z^2$, $Z^3$, which now determine the bracket-generating set of vector fields $S^1$, $S^2$, $S^3$ compeletely. 
The Lie algebra of split $\frak{g}_2$ that arises in this way can be viewed also as the symmetry algebra of the $(2,3,5)$-distribution annihilated by the 1-forms $\{\tilde \theta_1, \tilde \theta_2,\tilde \theta_3\}$ where
\begin{align*}
\tilde \theta_1&=\der \tilde a_1+\left(\tilde a_3+\frac{1}{2}\tilde a_4 \tilde a_5\right)\der \tilde a_4,\\
\tilde \theta_2&=\der \tilde a_2+\left(\tilde a_3-\frac{1}{2}\tilde a_4 \tilde a_5\right)\der \tilde a_5,\\
\tilde \theta_3&=\der \tilde a_3+\frac{1}{2}\tilde a_4 \der \tilde a_5-\frac{1}{2}\tilde a_5 \der \tilde a_4,
\end{align*}
under the transformation
\begin{align*}
(\tilde a_1,\tilde a_2,\tilde a_3,\tilde a_4,\tilde a_5)=(a_1+a_3 a_4,a_2+a_3 a_5,-a_3,a_4,a_5).
\end{align*}
The vector fields in the symmetry algebra are precisely the ones for which Lie derivative of $\tilde \theta_1$, $\tilde \theta_2$, $\tilde \theta_3$ with respect to these vector fields are in the span of $\{\tilde \theta_1, \tilde \theta_2,\tilde \theta_3\}$. See \cite{tw14} for more explanation. 

To illustrate the procedure for writing down the Lie algebra of split $\frak{g}_2$, let us look at the example of the Hilbert-Cartan distribution. To map the Hilbert-Cartan distribution given by the annihilator of the 1-forms
\[
\der y-p \der x,\quad \der p-q \der x,\quad \der z-q^2\der x,
\]
into the flat Cartan distribution, we take
\begin{align*}
a_1=2z+2q^2x-4pq,\quad a_2=2y,\quad a_3=2p-qx,\quad a_4=2q,\quad a_5=-x.
\end{align*}
This gives
\begin{align*}
\theta_1=2(\der z-q^2 \der x)-4q(\der p-q \der x),\quad \theta_2=2(\der y-p \der x), \quad \theta_3=2(\der p-q \der x).
\end{align*}
We determine
\begin{align*}
c_1&=12z-32pq+12q^2 x,\quad c_2=12y+4p x-4q^2x,\\
c_3&=4p-2qx,\quad c_4=-2q,\quad c_5=-x.
\end{align*}
Finding the basis of vector fields $\partial_{c_1}$, $\partial_{c_2}$, $\partial_{c_3}$, $\partial_{c_4}$, $\partial_{c_5}$, we obtain
\begin{align*}
Z^1&=\frac{1}{4}(\partial_p-x\partial_y+4q\partial_z),\\
Z^2&=\frac{x}{4}\left(\partial_p-x\partial_y+4q\partial_z\right)-\frac{1}{2}\partial_q,\\
Z^3&=-(\partial_x+p\partial_y+q^2\partial_z+q\partial_p)-\frac{q}{2}(\partial_p-x\partial_y+4q\partial_z),
\end{align*}
so that
\begin{align}
S^1&=-\frac{1}{2}\partial_q,\nonumber\\ \label{hcf1}
S^2&=-(\partial_x+p \partial_y+q^2\partial_z+q\partial_p),\\
S^3&=\frac{1}{2}(12z-32pq+12 q^2x)\partial_q-(2p-qx)^2(\partial_p-x\partial_y+4q\partial_z)\nonumber\\
&\quad{}-(12y+4 px-4 qx^2)(\partial_x+p \partial_y+q^2\partial_z+q\partial_p).\nonumber
\end{align}
The vector fields $S^1$ and $S^2$ are in the span of the Hilbert-Cartan distribution and together with $S^3$ they pairwise bracket-generate to form a split $\frak{g}_2$ Lie algebra. 
In this paper we compute the bracket-generating set of vector fields for the maximally symmetric distribution determined by solutions of the $k=\frac{2}{3}$ and $\frac{3}{2}$ generalised Chazy equation as discussed in \cite{r16} and \cite{r19}. Together with \cite{r21b} this establishes the equivalences of the maximally symmetric rolling distribution, maximally symmetric generalised Chazy distribution and the Hilbert-Cartan distribution to one another. 

\section{Generalised Chazy equation and the relationship to Lam\'e equation}

It was shown in an earlier work \cite{r16} that the generalised Chazy equation with parameters $k=\frac{3}{2}$ and $k=\frac{2}{3}$ occur in $(2,3,5)$-distributions with maximal $G_2$ symmetry. In the subsequent sections we shall make the link explicit by showing the coordinate changes that map such $(2,3,5)$-distributions to the flat Cartan distribution. The Lie algebra of vector fields is then the push-forward of the $\frak{g}_2$ vector fields associated to the flat Cartan distribution via the inverse of this local coordinate diffeomorphism. This allows us to write the vector fields that bracket generate $\frak{g}_2$ in terms of the generalised Chazy distribution. 

The generalised Chazy equation or Chazy XII equation with parameter $k$ is given by
\begin{align}\label{chazyg}
y'''-2yy''+3y'^2-\frac{4}{36-k^2}(6y'-y^2)^2=0
 \end{align} 
and Chazy's equation or the Chazy III equation
\begin{align}\label{chazy}
y'''-2yy''+3y'^2=0
 \end{align} 
is obtained in the limit as $k$ tends to infinity. We may refer to equation (\ref{chazy}) as the generalised Chazy equation with parameter $k=\infty$.
Here $'$ denotes differentiation with respect to $x$. The generalised Chazy equation was introduced in \cite{chazy1}, \cite{chazy2} and studied more recently in \cite{co96}, \cite{acht} and \cite{ach}. For our purposes here, we use the Lam\'e parameterisation of solutions of the generalised Chazy equation, following \cite{co96} and \cite{r18}. For the generalised Chazy equation written as a first order system, see \cite{r16} and \cite{r18b}.

The generalised Chazy equation can also be solved in terms of hypergeometric functions \cite{r16} \cite{r16b}. The solutions of the generalised Chazy equation with parameters $k=2$, $k=3$, $k=\frac{2}{3}$ and $k=\frac{3}{2}$ given by hypergeometric functions have already been explored in \cite{r16b} and in the examples later on, it can be checked that the results hold as well for the hypergeometric parametrisations. 

Taking $P=y$, the generalised Chazy equation (\ref{chazyg}) with parameter $k$ is equivalent to the following non-linear first order system of differential equations
\begin{align}\label{ndek0}
\frac{\der P}{\der x}&=\frac{1}{6}(P^2-Q),\nonumber\\
\frac{\der Q}{\der x}&=\frac{2}{3}(PQ-R),\\
\frac{\der R}{\der x}&=PR+\frac{k^2}{36-k^2}Q^2.\nonumber
\end{align}
For brevity call $\al=\frac{k^2}{k^2-36}$. We will be interested in the case where $k=\frac{3}{2}$, which gives $\al=-\frac{1}{15}$, and $k=\frac{2}{3}$, which gives $\al=-\frac{1}{80}$.

We shall discuss a method introduced in \cite{co96} of integrating the solutions to obtain the Lam\'e form of the solutions. For further details about Lam\'e equations see \cite{ww}. See also Section 2 of \cite{r18}. Let $\De=e^{2\int P \der x}$. Integrating the second and third equations of the system (\ref{ndek0}) gives us
\begin{align*}
Q&= \De^{\frac{1}{3}}\left(\int-\frac{2}{3}R\De^{-\frac{1}{3}}\der x+c_1\right),\\
R&=\De^{\frac{1}{2}}\left(\int-\al Q^2\De^{-\frac{1}{2}}\der x+c_2\right).
\end{align*}
Let us denote $\mu=\frac{Q}{\De^{\frac{1}{3}}}-c_1$, $\nu=\frac{R}{\De^{\frac{1}{2}}}-c_2$. We obtain the differential equations
\begin{align*}
\mu'&=-\frac{2}{3}R\De^{-\frac{1}{3}}=-\frac{2}{3}(\nu+c_2)\De^{\frac{1}{6}},\\
\nu'&=-\al Q^2\De^{-\frac{1}{2}}=-\al(\mu+c_1)^2\De^{\frac{1}{6}}.
\end{align*}
We now introduce the new coordinate $\tilde z=\int\De^{\frac{1}{6}}\der x$, which gives $\der \tilde z=\De^{\frac{1}{6}} \der x$  or $\frac{\der}{\der \tilde z}=\De^{-\frac{1}{6}} \frac{\der}{\der x}$, 
so that the equations become
\begin{align*}
\mu_{\tilde z}&=-\frac{2}{3}(\nu+c_2),\\
\nu_{\tilde z}&=-\al(\mu+c_1)^2,
\end{align*}
or that
\begin{align*}
\mu_{\tilde z \tilde z}=\frac{2}{3}\al(\mu+c_1)^2.
\end{align*}
Hence some multiple of $\mu$ satisfies the Weierstrass differential equation with $g_2=0$, and by translation we shall set $c_1$, or equivalently $g_1=0$.
Specifically, $\mu=\frac{9}{\al}\wp$ where the Weierstrass $\wp$ function satisfies the Weierstrass differential equation
\begin{equation*}
(\wp_{\tilde z})^2=4\wp^3-g_3. 
\end{equation*}
We also obtain $\nu=-\frac{27}{2\al}\wp_{\tilde z}-c_2$. The coordinate $\tilde z$ here is different from the coordinate labelled $z$ in $M_{xypqz}$ in the theory of $(2,3,5)$-distributions, and this complex coordinate will show up again subsequently only in the context of discussing solutions to the spin $\frac{3}{2}$ and spin $4$ Lam\'e equations in Sections \ref{spin32} and \ref{spin4} respectively. 

Now take $\Phi=\De^m$. Using $\De'=2P\De$, we find
\begin{align*}
\Phi_{\tilde z}&=m \De^{-\frac{1}{6}}(2 P \De) \De^{m-1}=2m P \De^{m-\frac{1}{6}}, \\
\Phi_{\tilde z \tilde z}&=\De^{-\frac{1}{6}}\left(2m P' \De^{m-\frac{1}{6}}+2m (m-\frac{1}{6})P (2P\De)\De^{m-\frac{7}{6}}\right)\\
&=\frac{1}{3}m (P^2-Q)\De^{m-\frac{1}{3}}+4m (m-\frac{1}{6})P^2 \De^{m-\frac{1}{3}}\\
&=\left(\frac{1}{3}m +4m (m-\frac{1}{6})\right)P^2\De^{m-\frac{1}{3}}-\frac{1}{3}Qm\De^{m-\frac{1}{3}}.
\end{align*}
We therefore eliminate terms involving $P^2$ if we take $\frac{1}{3}m +4m (m-\frac{1}{6})=4m^2-\frac{1}{3}m=0$, or $m=\frac{1}{12}$. 
This gives the Lam\'e equation for $\Phi$ to satisfy. We obtain 
\begin{align*}
\Phi_{\tilde z \tilde z}=-\frac{1}{36}Q\De^{-\frac{1}{4}}=-\frac{1}{36}\frac{Q}{\De^{\frac{1}{3}}}\Phi=-\frac{1}{36}\mu\Phi=-\frac{1}{4\al}\wp\Phi.
\end{align*}

We have the following:
\begin{theorem}\label{lame1}
Given a solution to the Lam\'e equation
\begin{align*}
\Phi_{\tilde z \tilde z}+\frac{(k+6)(k-6)}{4k^2}\wp\Phi=0
\end{align*}
where $\wp=\wp(\tilde z,0,g_3)$, the solution to the first order system (\ref{ndek0}) associated to the generalised Chazy equation with parameter $k$ can be parametrised by
\[
(P,Q,R)=\left(6\Phi_{\tilde z}\Phi,\frac{9(k^2-36)}{k^2}\wp \Phi^4,-\frac{27}{2}\frac{k^2-36}{k^2}\wp_{\tilde z}\Phi^6\right).
\]
\end{theorem}
 
When the Lam\'e equation is written in the form 
\begin{align*}
\Phi_{\tilde z \tilde z}-n(n+1)\wp\Phi=0,
\end{align*}
the value $n$ is called the spin of the Lam\'e equation. The Lam\'e equation that we consider in this article is called equianharmonic since $g_2=0$ and the equation also has zero accessory parameter. The correspondence between the Chazy parameter $k$ and spin $n$ is as follows:
\begin{equation*}
n=-\frac{1}{2}-\frac{3}{k} \text{~or~} n=-\frac{1}{2}+\frac{3}{k}.
\end{equation*}
Restricting to postive values of $k$ and values of $n \geq -\frac{1}{2}$, we have the following correspondence between the spin of the Lam\'e equation and the generalised Chazy parameter.
For $n=-\frac{1}{2}$, we obtain the critical Lam\'e equation. This corresponds to the case where $k=\infty$, or the Chazy equation (see also Section 2 of \cite{r18}). 
For half-integer spin, we have the following values for $k$.
When $n=\frac{1}{2}$, we have $k=3$. 
When $n=\frac{3}{2}$, we have $k=\frac{3}{2}$. For these values of $n$, the corresponding Lam\'e equation was studied in part two/Section 3 of \cite{r18}.
For integer spin, we have the following values for $k$.
When the spin $n=1$, we have $k=2$.
When $n=2$, we have $k=\frac{6}{5}$.
When $n=3$, we have $k=\frac{6}{7}$.
When $n=4$, we have $k=\frac{6}{9}=\frac{2}{3}$.
In this article we are interested in the cases where $k=\frac{2}{3}$, $\frac{3}{2}$, $2$ and $3$ which corresponds to spin $4$, $\frac{3}{2}$, $1$ and $\frac{1}{2}$ respectively. We give the solutions to the spin $\frac{3}{2}$ Lam\'e equation in Section \ref{spin32} and solutions to the spin $4$ Lam\'e equation in Section \ref{spin4}.

\section{Local equivalence of the maximally symmetric $k=\frac{3}{2}$ generalised Chazy distribution to flat Cartan distribution}

We consider the rank 2 distribution $\cD_{\varphi(x,q)}$ on $M_{xyzpq}$ associated to the underdetermined differential equation $z'=\varphi(x,y'')$ where $\varphi(x,y'')=\frac{(y'')^2}{H''(x)}$ and $H''(x)$ is a non-zero function of $x$. This is to say that the distribution $\cD_{\varphi(x,q)}$ is annihilated by the three 1-forms
\begin{align}\label{ch321forms}
\om_1&=\der y-p \der x,\nonumber\\
\om_2&=\der p-q \der x,\\
\om_3&=\der z-\frac{q^2}{H''(x)} \der x.\nonumber
\end{align}
In \cite{r19}, it was shown that the three 1-forms can be completed to a coframing $(\theta_1, \theta_2, \theta_3, \theta_4, \theta_5)$ on $M_{xyzpq}$ and this coframing gives a metric 
\[
g=2\theta_1\theta_5-2\theta_2\theta_4+\frac{4}{3}\theta_3\theta_3 
\]
in Nurowski's conformal class \cite{conf}. The metric $g$ is conformally flat, i.e. the metric $g$ has vanishing Weyl tensor if and only if $H(x)$ is a solution to the 6th-order nonlinear differential equation 
\begin{align}\label{noth}
10H''^3H^{(6)}-70H''^2H'''H^{(5)}&-49H''^2H''''^2+280H''H'''^2H''''-175H'''^4=0.
\end{align}
This equation is called Noth's equation \cite{annur}. Making the substitution $H''=e^{\frac{2}{3}\int P \der x}$, we obtain the $k=\frac{3}{2}$ generalised Chazy equation for $P(x)$. In this case the distribution of the form $\cD_{\varphi(x,q)}$ is maximally symmetric and in \cite{r19} the explicit form of the metric given by the distribution $\cD_{\varphi(x,q)}$ is written as well. It was also found in \cite{r19} that rescaling the metric $\tilde g=2^{-\frac{2}{3}}H''^{\frac{2}{3}}g$ further by a conformal factor $\Om$, we obtain a Ricci-flat representative in Nurowski's conformal class. That is to say ${\rm Ric}(\Om^2\tilde g)=0$. We find that $\Om^2 \tilde g$ is Ricci-flat when $\Om$ satisfies the second-order differential equation
\begin{equation}\label{rfa}
\Om''\Om-2\Om'^2-\frac{2}{3}P\Om\Om'-\frac{1}{18}P^2\Om^2-\frac{1}{30}Q\Om^2=0.
\end{equation}
We make the substitution $\Om=\frac{1}{\rho} e^{-\frac{1}{3}\int P\der x}$ to obtain
\begin{equation}\label{sode1}
\rho''-\frac{1}{45}Q\rho=0,
\end{equation}
where $Q=P^2-6P'$ and the solution $\rho(x)$ is given in Theorem 3.1 of \cite{r19}, involving both the solutions of the $k=3$ and $k=\frac{3}{2}$ generalised Chazy equation.

When
\[
\varphi(x,q)=\frac{q^2}{H''(x)},
\]
we can map the 1-forms given in (\ref{ch321forms}) into the 1-forms (\ref{t1}), (\ref{t2}) and (\ref{t3}) associated to the flat Cartan distribution as follows.
We define
\begin{align*}
\chi=\int \Om^2 H'' \der x. 
\end{align*}
Let us take
\begin{align*}
K&=\frac{1}{\Om H''},\\
L&=\frac{1}{\Om^2 H''^2}\left(\Om H'''+4\Om'H''\right),\\
M&=\frac{1}{\Om^2 H''^3}\left(\Om (H''H^{(4)}-H'''^2)+\Om'H''H'''+3\Om''H''^2\right).
\end{align*}
It can be checked that we have
\begin{align*}
\chi'&=\Om^2H'',\\
K'&=3\frac{\Om'}{\Om}K-L,\\
(L\chi')'&=\chi'M+\Om''.
\end{align*}
Define
\begin{align*}
a_1&=-\frac{\Om \chi'}{16}y,\\
a_2&=z+\bigg(M^2\chi-5\chi' L M+12\Om' M+\frac{1}{2\Om}\chi'^2 L^3+\frac{4\Om'}{\Om}\chi'L^2+(4\Om \Om''-34\Om'^2)\frac{1}{\Om}L\\
&\quad{}+\frac{1}{2}\frac{H^{(5)}}{H''^2}+\frac{2}{\Om^2}K(\Om'''\Om^2-10\Om\Om'\Om''+30\Om'^3)\bigg)y^2\\
&\quad{}+(L^2\chi-L\Om)p^2+K^2 \chi q^2+(2LM\chi-3\Om M-\Om' L-\Om'' K)yp\\
&\quad{} +2K(M\chi-\chi'L+\Om')yq+2K(L\chi-\Om)pq,\\
a_3&=(-\frac{1}{4}M\chi+\frac{1}{2}\chi' L-\frac{1}{2}\Om')y+(\frac{1}{2}\Om-\frac{1}{4}\chi L)p-\frac{1}{4}\chi Kq,\\
a_4&=\frac{1}{8}\chi,\\
a_5&=4My+4Lp+4Kq.
\end{align*}

We find for this set of functions, 
\begin{align*}
\theta_1&=-\frac{1}{16}\Om \chi' \om_1,\\
\theta_2&=\bigg(\big(\frac{H^{(5)}}{H''^2}+\frac{\chi'^2}{\Om}L^3-8L(M-\frac{\Om'}{\Om}L)\chi'+22M \Om'-4K(10\frac{\Om'\Om''}{\Om}-30\frac{\Om'^3}{\Om^2}-\Om''')\\
&\quad{}-4(17\frac{\Om'^2}{\Om}-2\Om'')L\big)y-(\Om'' K+\Om' L+\Om M)p+2K(\Om'-\chi' L)q\bigg)\om_1\\
&\quad{}+((2L^2\chi'-\Om'' K-3\Om M-3\Om' L)y-2\Om K q)\om_2+ \om_3\\
&\quad{}-\frac{yq}{\Om^2H''^3} b_1\der x-\frac{y^2}{2\Om^4H''^5}b_2\der x,\\
\theta_3&=\frac{1}{2}(\chi' L-\Om')\om_1+\frac{\Om}{2}\om_2,
\end{align*}
where
\begin{align*}
b_1=10H''^2\Om''\Om-20H''^2\Om'^2-10H''H'''\Om\Om'+3H''H''''\Om^2-5H'''^2\Om^2
\end{align*}
and
\begin{align*}
b_2&=(-H''^3H^{(6)}+8H''^2H'''H^{(5)}+8H''^2H''''^2-43H''H'''^2H''''+30H'''^4)\Om^4\\
&\quad{}+(-4H''^2\Om''''+14H''H'''\Om'''+26H''H''''\Om''-50H'''^2\Om'')H''^2\Om^3\\
&\quad{}+(4H''^2H^{(5)}-40H''H'''H''''+50H'''^3)H''\Om'\Om^3+38H''^4\Om''^2\Om^2\\
&\quad{}+(24H''^2\Om'\Om'''-98H''H'''\Om'\Om''-44H''H''''\Om'^2+100H'''^2\Om'^2)H''^2\Om^2\\
&\quad{}-40(4H''\Om''-3H'''\Om')H''^3\Om'^2\Om+120H''^4\Om'^4.
\end{align*}
The equation $b_1=0$ is precisely the equation for Ricci-flatness (\ref{rfa}) upon substituting $H(x)=\iint e^{\frac{2}{3}\int P(x) \der x}\der x \der x$. Solving $b_1=0$ for $\Om''$ and substituting it into the equation $b_2=0$ gives Noth's equation
\begin{align*}
H^{(6)}H''^3-7H''^2H'''H^{(5)}-\frac{49}{10}H''^2H''''^2+28H''H'''^2H''''-\frac{35}{2}H'''^4=0. 
\end{align*}
It follows that $\theta_1$, $\theta_2$ and $\theta_3$ are in the span of $\om_1$, $\om_2$, $\om_3$ if and only if both (\ref{noth}) and (\ref{rfa}) hold. We now assume that it is the case where $b_1=0$ and $b_2=0$.
The equation for $\Om''$ is equivalent to the equation
\begin{align*}
L'=-\frac{\Om'}{\Om}L+\frac{L^2}{K}-2M.
\end{align*}
Consequently, we also find that
\begin{align*}
M'&=K\left(\frac{5}{2}\frac{L^3}{K^3}-14\frac{\Om'}{\Om}\frac{L^2}{K^2}+10\frac{\Om'^2}{\Om^2}\frac{L}{K}+12\frac{\Om'^3}{\Om^3}+19\frac{\Om'}{\Om}\frac{M}{K}-4\frac{L}{K}\frac{M}{K}+\frac{1}{10}\frac{H^{(5)}}{H''}\right).
\end{align*}
To derive the vector fields that generate $\frak{g}_2$, we compute
\begin{align*}
(c_1,c_2,c_3,c_4,c_5)=\left(6a_1-2a_3a_4+a_4^2a_5,6a_2-2a_3a_5-a_4a_5^2,2a_3,-a_4,a_5\right).
\end{align*}

We have
\begin{align*}
\Om''K=2L^2\chi'-3\Om M-3\Om' L.
\end{align*}
Substituting this and its differential consequence into $(c_1,c_2,c_3,c_4,c_5)$ and using the formulas for $K'$, $L'$ and $M'$ gives
\begin{align*}
c_1&=\frac{\chi^2}{8}(My+Lp+Kq)-\frac{\chi}{8}((\chi'L-\Om')y+\Om p)-\frac{3\chi'}{8} \Om y,\\
c_2&=6z+6\chi(My+Lp+Kq)^2-16(My+Lp+Kq)\big((\chi'L-\Om')y+\Om p\big)\\
&\quad{}-3\left(\frac{1}{5}\frac{H^{(5)}}{H''^2}-10\frac{H'''}{H''}\Om M+K \Om \left(5\frac{L^3}{K^3}-28\frac{L^2}{K^2}\frac{\Om'}{\Om}+20\frac{L}{K}\left(\frac{\Om'}{\Om}\right)^2+24\left(\frac{\Om'}{\Om}\right)^3\right)\right)y^2\\
&\quad{}+6 \Om Lp^2+12\Om M yp,\\
c_3&=-\frac{\chi}{2}(My+Lp+Kq)+(\chi'L-\Om')y+\Om p,\\
c_4&=-\frac{1}{8}\chi,\\
c_5&=4(My+Lp+Kq).
\end{align*}
Assuming $b_1=0$ and $b_2=0$, the 1-forms given by $\Theta_1$, $\Theta_2$ and $\Theta_3$ are in the span of $\om_1$, $\om_2$ and $\om_3$.
This set of functions allows us to write down the split $\frak{g}_2$ Lie algebra of vector fields readily following Section 2. 

We now pass to the generalised Chazy equation by making the substitution 
\begin{align*}
\Om&=\frac{1}{\rho}\exp(-\frac{1}{3}\int P(x) \der x)\hspace{12pt} \mbox{and} \hspace{12pt} H(x)=\iint \exp(\frac{2}{3}\int P(x)\der x)\der x \der x,
\end{align*}
where $P(x)$ satisfies the $k=\frac{3}{2}$ generalised Chazy equation and we make use of the first-order system
\begin{align*}
P'&=\frac{1}{6}(P^2-Q),\\
Q'&=\frac{2}{3}(PQ-R),\\
R'&=PR+\frac{1}{15}Q^2,
\end{align*}
along with the fact that $\rho(x)$ satisfies the second-order differential equation (\ref{sode1})
\begin{align*}
\rho''-\frac{1}{45}Q\rho=0.
\end{align*}
Let us take
\[
\xi=\frac{\rho'}{\rho},
\]
so that
\[
\rho=e^{\int \xi \der x}.
\]
We find
\[
\chi=\int \frac{1}{\rho^2} \der x
\]
and define
\begin{align*}
J=My+Lp+Kq. 
\end{align*}
We obtain
\begin{align*}
J=\frac{2}{3}\rho e^{-\frac{1}{3}\int P \der x}\left(9y \xi^2+2(Py-3p)\xi-(\frac{Q}{60}-\frac{P^2}{12})y-Pp+\frac{3}{2}q\right),
\end{align*}
which gives
\begin{align*}
c_1&=\frac{\chi^2}{8}J+\frac{\chi}{8\rho}e^{-\frac{1}{3}\int P \der x}\left(3y\xi+\frac{P}{3}y-p\right)-\frac{3}{8\rho^3}e^{-\frac{1}{3}\int P \der x}y,\\
c_2&=6z+6\chi J^2+\frac{16}{\rho} J e^{-\frac{1}{3}\int P \der x}(3 y \xi+\frac{P}{3}y-p)+e^{-\frac{2}{3}\int P \der x}\bigg(-4(6\xi+P) p^2\\
&\quad{}+2(36 \xi^2+8 P \xi-\frac{1}{15}Q+\frac{1}{3}P^2) yp-(72\xi^3+24P\xi^2+\frac{8}{3}P^2\xi+\frac{1}{9}(P^3+\frac{1}{5}PQ+\frac{2}{5}R))y^2\bigg),\\
c_3&=-\frac{\chi}{2}J-\frac{1}{\rho}e^{-\frac{1}{3}\int P \der x}\left(3y\xi+\frac{P}{3}y-p\right),\\
c_4&=-\frac{1}{8}\chi,\\
c_5&=4J.
\end{align*}

Let us take
\begin{align}\label{x1a}
X^1&=\partial_q,\\ \label{x2a}
X^2&=\partial_x+p\partial_y+q\partial_p+q^2e^{-\frac{2}{3}\int P \der x}\partial_z,
\end{align}
to be the vector fields that span the $k=\frac{3}{2}$ generalised Chazy distribution, which are also annihilated by the 1-forms in (\ref{ch321forms}). We also take
\begin{align*}
X^3&=4\rho^3\chi e^{-\frac{1}{3}\int P \der x}\bigg(3y\xi^3+3(\frac{P}{3}y-p)\xi^2-((\frac{Q}{60}-\frac{P^2}{12})y+Pp)\xi\\
&\quad{}-\frac{1}{270}(PQ+R)y-\frac{1}{12}(P^2+\frac{Q}{15})p\bigg)\partial_z-\rho^3\chi e^{\frac{1}{3}\int P \der x}(\partial_y+\frac{1}{3}(P+9\xi)\partial_p)\\
&\quad{}+4\rho e^{-\frac{1}{3}\int P \der x}\left(3y\xi^2+2(\frac{P}{3}y-p)\xi-\frac{1}{3}((\frac{Q}{60}-\frac{P^2}{12})y+pP-3q)\right)\partial_z\\
&\quad{}+\rho e^{\frac{1}{3}\int P \der x}\partial_p.
\end{align*}
This gives 
\begin{align*}
Z^1&=-\frac{\rho e^{\frac{1}{3}\int P dx}}{6}(\big(36\xi^2+12P\xi+P^2+\frac{Q}{15}\big)\rho^2 \chi-4 P-24\xi)X^1+X^3,\\
Z^2&=\ell_1X^1-8\rho^2X^2-\frac{8}{3}\rho e^{-\frac{1}{3}\int P \der x}(9y\xi^2+6(\frac{P}{3}y-p)\xi-(\frac{Q-5P^2}{60}y+Pp-\frac{3}{2}q))X^3,\\
Z^3&=\ell_2X^1-\frac{1}{8}\chi X^3,
\end{align*}
where the functions $\ell_1$ and $\ell_2$ are given by
\begin{align*}
\ell_1&=\frac{4\rho^4}{9}\big(36\xi^2+12P\xi+P^2+\frac{Q}{15}\big)\big(9y\xi^2+6(\frac{P}{3}y-p)\xi-(\frac{Q-5P^2}{60})y-Pp+\frac{3}{2}q\big) \chi \\
&\quad{}+\rho^2\big(-144y\xi^3+16(-\frac{10}{3}Py+7p)\xi^2+\frac{4}{3}((\frac{1}{3}Q-\frac{13}{3}P^2)y+28Pp-30q)\xi\\
&\quad{}+\frac{4}{9}((\frac{2}{15}R-\frac{1}{3}P^3+\frac{1}{5}QP)y+(7P^2+\frac{1}{5}Q)p-24Pq)\big),\\
\ell_2&=\frac{e^{\frac{1}{3}\int P \der x}}{48\rho}((36\xi^2+12P\xi+P^2+\frac{Q}{15})\rho^4\chi^2-4(P+6\xi)\rho^2\chi+12).
\end{align*}
Our preference for ordering the terms is given by $\rho^2 \chi$, $\xi$, $y$, $p$, $q$.
Finally, we obtain
\begin{align}\label{s1a}
S^1&=-\rho^2\bigg(48y\xi^3+48(\frac{P}{3}y-p)\xi^2-16((\frac{Q-5P^2}{60})y+pP-\frac{3}{2}q)\xi\\
&\quad{}-\frac{4}{3}(\frac{2}{45}(PQ+R)y+(P^2+\frac{Q}{15})p-6Pq)\bigg)X^1-8\rho^2X^2,\nonumber \\\label{s2a}
S^2&=\frac{e^{\frac{1}{3}\int P \der x}}{4\rho}X^1.
\end{align}
These two vector fields lie in the span of $X^1$ and $X^2$ and together with 
\begin{align}\label{s3a}
S^3&=\ell_3X^1+\frac{2e^{-\frac{1}{3}\int P \der x}}{3\rho}\bigg(\rho^4\chi^2(9y \xi^2+2(Py-3p)\xi-(\frac{Q}{60}-\frac{P^2}{12})y-Pp+\frac{3}{2}q)\\
&\quad{}+\frac{3}{2}\rho^2 \chi (3y \xi+\frac{P}{3}y-p)-\frac{9}{2}y\bigg)X^2\nonumber\\
&\quad{}-\frac{e^{-\frac{2}{3}\int P \der x}}{9\rho^2}(\rho^2\chi (9y \xi^2+2(Py-3p)\xi-(\frac{Q}{60}-\frac{P^2}{12})y-Pp+\frac{3}{2}q)+9 y\xi+Py-3p)^2X^3, \nonumber
\end{align} 
where
\begin{align*}
\ell_3&=\frac{\rho e^{\frac{1}{3}\int P dx}}{6}\big((36\xi^2+12P\xi+P^2+\frac{Q}{15})\rho^2 \chi-4 P-24\xi\big)c_3^2\\
&\quad{}+\frac{e^{\frac{1}{3}\int P \der x}}{4\rho}c_2+\rho^2\bigg(48\xi^3y+48(\frac{P}{3}y-p)\xi^2-16((\frac{Q-5P^2}{60})y+pP-\frac{3}{2}q)\xi\\
&\quad{}-\frac{4}{3}(\frac{2}{45}(PQ+R)y+(P^2+\frac{Q}{15})p-6Pq)\bigg)c_1,
\end{align*}
they bracket-generate the Lie algebra of split $\frak{g}_2$. This proves the following result. 

\begin{theorem}\label{g2a}
For the maximally symmetric distribution spanned by the vector fields $X^1$ in (\ref{x1a}) and $X^2$ in (\ref{x2a}) where $P(x)$ satisfies the generalised Chazy equation with parameter $k=\frac{3}{2}$ and annihilated by the 1-forms in (\ref{ch321forms}), the Lie algebra of split $\frak{g}_2$ is obtained from the pairwise bracket-generating set $\{S^1,S^2,S^3\}$ where $S^1$, $S^2$ and $S^3$ are given in (\ref{s1a}), (\ref{s2a}), (\ref{s3a}).
\end{theorem}

\section{Examples: Hilbert-Cartan equation and spin $\frac{3}{2}$ Lam\'e equation}\label{spin32}
Specialising to the Hilbert-Cartan distribution obtained when $(P,Q,R)=(0,0,0)$, we see that there are non-constant conformal rescalings of the flat metric to itself given by $\rho=\al x+\beta$, where $\al$ and $\beta$ are arbitrary constants both not simultaneously zero. In this case, we have
\begin{align}\label{hc3}
\chi=\int \frac{1}{(\al x+\beta)^2}\der x=-\frac{1}{\al(\al x+\beta)}+c \mbox{~and~} \xi=\frac{\al}{\al x+\beta}.
\end{align}
This gives
\begin{align*}
X^1&=\partial_q,\\
X^2&=\partial_x+p\partial_y+q\partial_p+q^2\partial_z,\\
X^3&=12\rho^3\chi \xi^2(y \xi-p)\partial_z-\rho^3\chi(\partial_y+3\xi \partial_p)+\rho\partial_p+4\rho(3y\xi^2-2p\xi+q)\partial_z.
\end{align*}

We obtain the following corollary:
\begin{corollary}
There is a three parameter family of vector fields $\{S^1,S^2,S^3\}$ associated to the Hilbert-Cartan distribution that bracket-generates the Lie algebra of split $\frak{g}_2$. They are given by 
\begin{align*}
S^1&=-24\rho^2(2\xi^2(y\xi-p)+q \xi)X^1-8\rho^2X^2,\\
S^2&=\frac{1}{4\rho}X^1,\\
S^3&=\ell_3X^1+\frac{2}{\rho}(\rho^4\chi^2(3y\xi^2-2p\xi+\frac{1}{2}q)-\frac{3}{2}y+\frac{1}{2}\rho^2\chi(3y\xi-p))X^2\\
&\quad{}-\frac{1}{\rho^2}(\rho^2\chi(3y\xi^2-2p\xi+\frac{1}{2}q)+3y\xi-p)^2X^3
\end{align*}
with
\begin{align*}
c_1&=\frac{\chi^2}{8}(6y\xi^2-4p\xi+q)\rho+\frac{\chi}{8\rho}(3y\xi-p)-\frac{3}{8\rho^3}y,\\
c_2&=6z+6\chi \rho^2(6y \xi^2-4p \xi+q)^2+16(6y\xi^2-4p \xi+q)(3y\xi-p)-24\xi(p^2-3\xi y p+3\xi^2 y^2),\\
c_3&=-\frac{\chi}{2}\rho(6y\xi^2-4p\xi+q)-\frac{1}{\rho}(3y\xi-p),\\
\ell_3&=\rho(6\xi^2\rho^2 \chi-4\xi)c_3^2+\frac{1}{4\rho}c_2+24\rho^2(2\xi^3y-2\xi^2p+q\xi)c_1,
\end{align*}
and where $\rho=\al x+\beta$, $\chi$ and $\xi$ are given as in (\ref{hc3}).
\end{corollary}
The vector fields (\ref{hcf1}) are recovered up to a constant multiple when we specialise to $\al=0$, $\beta=-\frac{1}{2}$, $c=0$, so that $\xi=0$.

Turning to the spin $\frac{3}{2}$ equianharmonic Lam\'e equation, recall that Halphen's solution to the spin $\frac{3}{2}$ Lam\'e equation (again with $g_2=0$)
\begin{align}\label{lame32}
\Phi_{\tilde z \tilde z}-\frac{15}{4}\wp(\tilde z,0,g_3)\Phi=0
\end{align}
is given as follows, see pages 379-380 of \cite{halphen} and also page 92 of \cite{r18}. 
Let $w=\wp\left(\frac{\tilde z}{2}\right)$. The solution to (\ref{lame32}) is given by 
\[
\Phi=Y(w)\wp'\left(\frac{\tilde z}{2}\right)^{-\frac{3}{2}}
\]
where
\[
Y(w)=\al w+\beta (w^3+\frac{g_3}{2}). 
\]
In this context $\wp'$ denotes differentiation with respect to the independent variable $\frac{\tilde z}{2}$.
Since
\[
\wp'\left(\frac{\tilde z}{2}\right)=\sqrt{4w^3-g_3},
\]
we therefore obtain
\begin{align*}
\Phi=\frac{\al w+\beta (w^3+g_3/2)}{(\wp'(\tilde z/2))^{\frac{3}{2}}}=\frac{\al w+\beta (w^3+g_3/2)}{(4w^3-g_3)^{\frac{3}{4}}}.
\end{align*}
From Theorem \ref{lame1} (see also \cite{r18}), we can parametrise the solutions of the $k=\frac{3}{2}$ Chazy equation by solutions of the spin $\frac{3}{2}$ Lam\'e equation as follows:
\[
(P,Q,R)=(6\Phi_{\tilde z}\Phi,-135\wp \Phi^4,\frac{405}{2}\wp_{\tilde z}\Phi^6),
\]
where 
\[
x=\int \frac{1}{\Phi^2}\der \tilde z. 
\]
We have $\der x=\frac{1}{\Phi^2}\der \tilde z$.
We obtain the parametrisation of $H''(x)$ by the new independent variable $\tilde z$ as follows:
\begin{align*}
H''(x)=\exp(-\frac{2}{3}\int P\der x)&=\exp(-\frac{2}{3}\int 6\frac{\Phi_{\tilde z}}{\Phi}\der \tilde z)\\
&=\exp(-4\ln \Phi)=\Phi^{-4}.
\end{align*}
This allows us to parametrise the 1-forms in (\ref{ch321forms}) as
\begin{align*}
\der y-\frac{p}{\Phi^2}\der \tilde z, \hspace{12pt} \der p-\frac{q}{\Phi^2}\der \tilde z, \hspace{12pt} \der z-\frac{q^2}{\Phi^6}\der \tilde z. 
\end{align*}
We actually want to work with the independent coordinate $w$, so we make a further change of coordinates, using 
\[
\der w=\frac{1}{2}\wp'(\tilde z/2)\der \tilde z
\]
or equivalently,
\[
\der \tilde z=\frac{2}{\wp'(\tilde z/2)}\der w.
\]
It follows that
\[
\der x=\frac{1}{\Phi^2}\der \tilde z=\frac{2}{Y(w)^2}(\wp')^2(\tilde z/2)\der w=\frac{2}{Y(w)^2}(4w^3-g_3)\der w.
\]
By making this change of coordinates, we obtain the Pfaffian system given by the 1-forms
\begin{align}\label{lame32-1forms}
\om_1&=\der y-\frac{2p}{Y^2}(4 w^3-g_3) \der w,\nonumber\\
\om_2&=\der p-\frac{2q}{Y^2} (4 w^3-g_3)\der w,\\
\om_3&= \der z-\frac{2q^2}{Y^6} (4 w^3-g_3)^4\der w,\nonumber
\end{align}
where $Y=\al w+\beta (w^3+\frac{g_3}{2})$.

To find the vector fields parametrised by solutions of the spin $\frac{3}{2}$ Lam\'e equation that bracket-generate the split Lie algebra of $\frak{g}_2$, we make the further substitution
\begin{align*}
\rho&=\frac{(\ga w+\de)\wp'(\tilde z/2)}{\al w+\beta(w^3+g_3/2)},\quad \chi=\int \frac{2}{(\ga w+\de)^2} \der w, \quad \xi=\frac{\Phi^2 \wp'(\tilde z/2)}{2}\frac{\rho_{w}}{\rho}, \\
\Phi&=\frac{\al w+\beta (w^3+g_3/2)}{\wp'(\tilde z/2)^{3/2}}, \quad e^{2\int P \der x}=\Phi^{12}.
\end{align*}
Observe that the solution of the spin $\frac{1}{2}$ Lam\'e equation is given by $\Phi=(\ga w+\de) \wp'(\tilde z/2)^{-\frac{1}{2}}$ (see \cite{r18}), so that $\rho$ can be viewed as a quotient of the spin $\frac{1}{2}$ solution over the spin $\frac{3}{2}$ solution, and $\chi$ can be viewed as the independent variable of the generalised Chazy equation with parameter $k=3$ (which corresponds to spin $\frac{1}{2}$). This agrees with the result of Theorem 3.1 in \cite{r19}. For simplicity, we now specialise to the case where $\ga=0$ and take the constant of integration $c=0$, so that $\chi=\frac{2}{\de^2}w$.

We have
\begin{align}\label{c1a}
c_1&=\frac{1}{4\de^3Y^3}(-2Y^3(\al w+\frac{3g_3}{8}\beta)y+Y^2(8w^3+g_3)wp+2w^2(4w^3-g_3)^2q),\\ \label{c2a}
c_2&=6z-12\beta(\al w+\frac{3g_3\beta}{16})y^2+\frac{12w^2(16w^3+5g_3)}{Y^2}p^2+12\frac{(4w^3-g_3)^4 w}{Y^6}q^2\\ \nonumber
&\quad{}+\frac{6w (\beta(8w^3+ g_3)-8\al w)}{Y}y p-\frac{4(4w^3-g_3)^2(2\al-3\beta w^2)}{Y^3}y q\\ \nonumber
&\quad{}+16\frac{(5w^3+g_3)(4w^3-g_3)^2}{Y^4}pq,\\ \label{c3a}
c_3&=\frac{1}{ \de Y^3}(\frac{1}{2}\al Y^3y-(2w^3+g_3)Y^2p-w(4w^3-g_3)^2q),\\ \nonumber
c_4&=-\frac{1}{4\de^2} w,\\ \nonumber
 c_5&=\frac{6\de}{Y^3}(\beta w Y^3 y+4w^2Y^2p+\frac{2}{3}(4w^3-g_3)^2q).
\end{align}
This gives
\begin{align*}
\Theta_1&=-\frac{1}{4\de^3}(\al w+\frac{3}{4} \beta g_3)\om_1+\frac{1}{4\de^3}\frac{w(4w^3-g_3)}{Y}\om_2,\\
\Theta_2&=\left(-6\beta(\al w+\frac{3}{4}\beta g_3)y-\frac{24w}{Y}(\al w+\frac{3}{4}\beta g_3)p-4\frac{(3\beta w^2+\al)(4w^3-g_3)^2}{Y^3}q\right)\om_1\\
&\quad{}+\frac{6(4w^3-g_3)}{Y}\left(\beta wy+\frac{4 w^2}{Y}p-\frac{4(4w^3-g_3)^2}{3Y^3}q\right)\om_2+6\om_3,\\
\Theta_3&=\frac{1}{2\de}(\al+3\beta w^2)\om_1+\frac{4w^3-g_3}{\de Y}\om_2,
\end{align*}
where $\om_1$, $\om_2$ and $\om_3$ are the 1-forms in (\ref{lame32-1forms}).
Using this set of functions as a basis, we find 
\begin{align*}
Z^1&=\partial_{c_3}+2c_5\partial_{c_2}-2c_4\partial_{c_1}\\
&=-\frac{2 \de w}{Y}\partial_y+\frac{4\de(4w^3-g_3)^2}{ Y^3}q\partial_z+\frac{\de(\beta (8w^3+g_3)+4 \al w)}{2(4w^3-g_3)}\partial_p-\frac{12 \de Y^3 w^2}{(4w^3-g_3)^3}\partial_q,
\end{align*}
and the corresponding formulas for $Z^2$ and $Z^3$, which we do not display here. 
We have the following corollary of Theorem \ref{g2a}.
\begin{corollary}\label{cora}
Consider the maximally symmmetric $(2,3,5)$-distribution given by the span of 
\begin{align*}
X^1&=\partial_q, \\
X^2&=\frac{Y^2}{2(4w^3-g_3)}\partial_w+p\partial_y+q\partial_p+\frac{(4w^3-g_3)^3}{Y^4}q^2\partial_z,
\end{align*}
where $Y=\al w+\beta (w^3+\frac{g_3}{2})$.
This is annihilated by the 1-forms $\{\om_1,\om_2,\om_3\}$ in (\ref{lame32-1forms}). Then the vector fields
\begin{align*}
S^1&=\bigg(\frac{6\beta \de^2 Y^3}{(4w^3-g_3)^2}y+\frac{24\de^2w Y^2}{(4w^3-g_3)^2}p+\frac{12\de^2(\al(8w^3+g_3)+9\beta g_3 w^2)}{(4w^3-g_3)Y}q\bigg)X^1\\
&\quad{}-\frac{8\de^2(4w^3-g_3)}{Y^2} X^2,\\
S^2&=\frac{Y^3}{4\de(4w^3-g_3)^2}X^1,\\
S^3&=-c_1S^1+c_2S^2-c_3^2Z^1,
\end{align*}
where 
\begin{align*}
Z^1&=-\frac{2 \de w}{Y}\partial_y+\frac{4\de(4w^3-g_3)^2}{Y^3}q\partial_z+\frac{\de(\beta (8w^3+g_3)+4 \al w)}{2(4w^3-g_3)}\partial_p-\frac{12 \de Y^3 w^2}{(4w^3-g_3)^3}\partial_q,
\end{align*}
and $c_1$, $c_2$ and $c_3$ are given in (\ref{c1a}), (\ref{c2a}) and (\ref{c3a}), pairwise bracket-generate the Lie algebra of split $\frak{g}_2$.
\end{corollary}

\section{Local equivalence of the maximally symmetric $k=\frac{2}{3}$ generalised Chazy distribution to flat Cartan distribution}

The function $H(x)$ is related to another function $F(\tilde x)$ by a Legendre transformation \cite{annur}, \cite{r16}. We say that $F(\tilde x)$ is the Legendre dual of $H(x)$ determined by the relation $H(x)+F(\tilde x)=x\tilde x$. This implies $\tilde x=H'(x)$ with $\der \tilde x=H''\der x$ and $H''=\frac{1}{F_{\tilde x\tilde x}}$. We can make use of this transformation to write $\der x=F_{\tilde x\tilde x}\der \tilde x$. The Legendre dual of the distribution $\cD_{\varphi(x,q)}$ is therefore given by the annihilator of the three 1-forms
\begin{align*}
\om_1&=\der y-p F_{\tilde x \tilde x}\der \tilde x,\nonumber\\
\om_2&=\der p-q F_{\tilde x \tilde x}\der \tilde x,\\
\om_3&=\der z-q^2 F_{\tilde x \tilde x}^2\der \tilde x\nonumber
\end{align*}
on the mixed jet space with local coordinates $(\tilde x, y, z, p, q)$. Relabelling $\tilde x$ with $x$, we have
\begin{align}
\om_1&=\der y-p F''\der x,\nonumber\\\label{ch231forms}
\om_2&=\der p-q F''\der x,\\
\om_3&=\der z-q^2 F''^2\der x.\nonumber
\end{align}
Here $F$ now becomes a function of $x$, and $'$ denotes differentiation with respect to $x$. These three 1-forms are completed to a coframing $(\theta_1, \theta_2, \theta_3, \theta_4, \theta_5)$ for a metric (\ref{metric}) in Nurowski's conformal class, as was done in \cite{r19}. The condition that the metric $g$ is conformally flat, i.e.\ the metric $g$ has vanishing Weyl tensor, occurs when $F(x)$ is a solution to the nonlinear differential equation 
\begin{align}\label{6thode}
10F''^3F^{(6)}-80F''^2F^{(3)}F^{(5)}&-51F''^2F''''^2+336F''F'''^2F''''-224F'''^4=0.
\end{align}
This equation appears in \cite{annur} and we can call it the dual of Noth's equation. 
If we replace $F''(x)=e^{\int \frac{1}{2}P(x)\der x}$, then equation (\ref{6thode}) is reduced to the generalised Chazy equation (\ref{chazyg}) for $y=P$ with parameter $k=\frac{2}{3}$. 

The metric $\tilde g=2^{\frac{1}{3}}(F'')^{-\frac{2}{3}}g$ can again be rescaled by a conformal factor to obtain a Ricci-flat representative in the conformal class \cite{r19}. We find that the Ricci tensor of $\Om^2 \tilde g$ is zero when $\Om$ satisfies
\begin{equation}\label{rf2b}
40\Om''\Om-80\Om'^2-6\Om^2P'+\Om^2P^2=0.
\end{equation}
If we make the substitution $\Om=\frac{1}{\eta}$, then we obtain the differential equation
\begin{equation}\label{sode2}
\eta''-\frac{1}{40}Q\eta=0
\end{equation}
where $Q=P^2-6P'$ and the solution $\eta$ was obtained in Theorem 3.2 of \cite{r19}. It involves both the solutions of the $k=2$ and $k=\frac{2}{3}$ generalised Chazy equation. 

We can map the 1-forms given by (\ref{ch231forms}) into the flat Cartan distribution as follows. Let us take
\begin{align*}
\chi&=\int \Om^2\der x,\\
K&=\frac{F''}{\Om},\\
L&=-\frac{1}{\Om}\left(\frac{F'''}{F''}-4\frac{\Om'}{\Om }\right),\\
M&=-\frac{1}{\Om^2 F''^3}\left(\Om (F''F^{(4)}-2F'''^2)+4\Om'F''F'''-3\Om''F''^2\right).
\end{align*}
It can be checked that we have
\begin{align*}
\chi'=\Om^2,\qquad K'=3\frac{\Om'}{\Om}K-\Om K L.
\end{align*}
Define the functions mapping into the 1-forms (\ref{t1}), (\ref{t2}), (\ref{t3}) annihilating the flat Cartan distribution by taking
\begin{align*}
a_1&=-\frac{\Om^2}{16K}y,\\
a_2&=z+\bigg(M^2\chi-6\frac{\Om'}{\Om}\frac{M}{K}-\frac{1}{2}\frac{F^{(5)}}{F''^3}-\frac{2}{K^2}\big(\Om L^3-\frac{10\Om'}{\Om}L^2+(\frac{5}{2} \frac{\Om''}{\Om^2}+27\frac{\Om'^2}{\Om^3})L\\
&\quad{}-\frac{1}{\Om^2}(\frac{\Om'''}{\Om}+5\frac{\Om'\Om''}{\Om^2}+22\frac{\Om'^3}{\Om^3})\big)\bigg)y^2+(L^2\chi-L\Om)p^2+K^2 \chi q^2\\
&\quad{}+(2LM\chi-3\Om M-2\frac{\Om'}{\Om}\frac{L}{K}+(4\frac{\Om'^2}{\Om^3}-\frac{\Om''}{\Om^2})\frac{1}{K})yp\\
&\quad{} +2K(M\chi-\Om \frac{L}{K}+\frac{\Om'}{\Om K})yq+2K(L\chi-\Om)pq,\\
a_3&=(-\frac{1}{4}M\chi+\frac{\Om}{2}\frac{L}{K}-\frac{1}{2}\frac{\Om'}{\Om K})y+(\frac{1}{2}\Om-\frac{1}{4}\chi L)p-\frac{1}{4}\chi Kq,\\
a_4&=\frac{1}{8}\chi,\\
a_5&=4My+4Lp+4Kq.
\end{align*}

We find for this set of functions, 
\begin{align*}
\theta_1&=-\frac{1}{16K}\Om^2\om_1,\\
\theta_2&=\bigg(\big(-\frac{F^{(5)}}{F''^3}+2\Om\frac{LM}{K}-14\frac{\Om'}{\Om}\frac{M}{K}-4\Om \frac{L^3}{K^2}+40\frac{\Om'}{\Om}\frac{L^2}{K^2}-(108\frac{\Om'^2}{\Om^3}+10\frac{\Om''}{\Om^2})\frac{L}{K^2}\\
&\quad{}+\frac{1}{K^2}(88\frac{\Om'^3}{\Om^5}+20\frac{\Om'\Om''}{\Om^4}+4\frac{\Om'''}{\Om^3})\big)y-(\Om M+2\frac{\Om'}{\Om}\frac{L}{K}-(4\frac{\Om'^2}{\Om^3}-\frac{\Om''}{\Om^2})\frac{1}{K})p\\
&\quad{}-2(\Om L-\frac{\Om'}{\Om})q\bigg)\om_1+\left(\big(2\Om \frac{L^2}{K}-4\frac{\Om'}{\Om}\frac{L}{K}-3\Om M+(4\frac{\Om'^2}{\Om^3}-\frac{\Om''}{\Om^2})\frac{1}{K}\big)y-2\Om K q\right)\om_2\\
&\quad{} + \om_3+\frac{y q}{\Om^2F''^2}b_1 \der x+\frac{y^2}{2\Om^4F''^6}b_2 \der x,\\
\theta_3&=\frac{1}{2}\left(\frac{\Om L}{K}-\frac{\Om'}{\Om K}\right)\om_1+\frac{\Om}{2}\om_2,
\end{align*}
where 
\begin{align*}
b_1=-10F''^2\Om''\Om+20F''^2\Om'^2-4F'''^2\Om^2+3F''F''''\Om^2
\end{align*}
and
\begin{align*}
b_2&=(-F''^3F^{(6)}+7F''^2F'''F^{(5)}+2F''^2F''''^2-20F''F'''^2F''''+12F'''^4)\Om^4\\
&\quad{}+(4F''^4\Om''''+10F''^3F''''\Om''-10F''^3F'''\Om'''-10F''^2F'''^2\Om'')\Om^3\\
&\quad{}+(20F''^2F'''^2\Om'^2-20F''^3F''''\Om'^2-24F''^4\Om'\Om'''+50F''^3F'''\Om'\Om''-38F''^4\Om''^2)\Om^2\\
&\quad{}-40F''^3F'''\Om\Om'^3+160F''^4\Om\Om'^2\Om''-120F''^4\Om'^4,
\end{align*}
so that $\theta_1$, $\theta_2$ and $\theta_3$ are in the span of $\om_1$, $\om_2$, $\om_3$ precisely when 
$b_1=0$ and $b_2=0$.
The equation $b_1=0$ is again the equation for Ricci-flatness (\ref{rf2b}).
Solving the equation $b_1=0$ for $\Om''$ and substituting it into $b_2=0$ gives the dual of Noth's equation
\begin{align*}
10F^{(6)}F''^3-80F''^2F'''F^{(5)}-51F''^2F''''^2+336F''F'''^2F''''-224F'''^4=0. 
\end{align*}

We have
\begin{align*}
L'&=\Om(KM-L^2)+3\frac{\Om'}{\Om}L+\frac{\Om''}{\Om^2}-4\frac{\Om'^2}{\Om^3},\\
M'&=6\Om L M-4\Om\frac{L^3}{K}-21\frac{\Om'}{\Om}M+24\frac{\Om'}{\Om}\frac{L^2}{K}-42\frac{\Om'^2}{\Om^3}\frac{L}{K}+\frac{28 \Om'^3}{\Om^5K}-\frac{F^{(5)}}{10 K^2 \Om^3}.
\end{align*}

Now assuming $b_1=0$ and $b_2=0$, to derive the vector fields that generate $\frak{g}_2$, we compute
\begin{align*}
(c_1,c_2,c_3,c_4,c_5)=\left(6a_1-2a_3a_4+a_4^2a_5,6a_2-2a_3a_5-a_4a_5^2,2a_3,-a_4,a_5\right)
\end{align*}
and find
\begin{align*}
c_1&=\frac{1}{8}(My+Lp+Kq)\chi^2+\frac{1}{8\Om K}(\Om' y-\Om^2(pK+Ly))\chi-\frac{3}{8}\frac{\Om^2}{K}y,\\
c_2&=6z+6(Kq+Lp+My)^2\chi-16(Kq+Lp+My)((\Om\frac{L}{K}-\frac{\Om'}{\Om K})y+\Om p)\\
&\quad{}+\bigg(\frac{3}{5}\frac{F^{(5)}}{\Om^3K^3}+\frac{12}{\Om^2 K^2}(2\Om^3L^3-12\Om^2L^2\frac{\Om'}{\Om}+21\Om L \frac{\Om'^2}{\Om^2}-14\frac{\Om'^3}{\Om^3})\\
&\quad{}-30(\Om L-4\frac{\Om'}{\Om})\frac{M}{K}\bigg)y^2+6\Om Lp^2+12\Om Myp,\\
c_3&=-\frac{1}{2}(My+Lp+Kq)\chi+(\Om \frac{L}{K}-\frac{\Om'}{\Om K})y+\Om p,\\
c_4&=-\frac{1}{8}\chi,\\
c_5&=4(My+Lp+Kq).
\end{align*}
The 1-forms given by $\Theta_1$, $\Theta_2$ and $\Theta_3$ are in the span of $\om_1$, $\om_2$ and $\om_3$ precisely when both $b_1=0$ and the dual of Noth's equation $b_2=0$ are satisfied. 

We now make the substitution 
\begin{align*}
\Om=\frac{1}{\eta}\hspace{12pt} \mbox{~and~} \hspace{12pt} F(x)=\iint \exp(\frac{1}{2}\int P(x)\der x)\der x \der x,
\end{align*}
where $(P,Q,R)$ satisfies
\begin{align*}
P'&=\frac{1}{6}(P^2-Q),\\
Q'&=\frac{2}{3}(PQ-R),\\
R'&=PR+\frac{1}{80}Q^2,
\end{align*}
and use equation (\ref{sode2})
\begin{align*}
\eta''-\frac{1}{40}Q\eta=0. 
\end{align*}
Let us again denote 
\[
\xi=\frac{\eta'}{\eta},
\]
so that
\[
\eta=e^{\int \xi \der x}.
\]
We get 
\[
\chi=\int \frac{1}{\eta^2} \der x
\]
and we obtain
\begin{align*}
J=My+Lp+Kq=\frac{\eta}{6}e^{-\frac{1}{2}\int P \der x}((6\xi+P)^2+\frac{Q}{20})y-\frac{\eta}{2}(8\xi+P)p+\eta e^{\frac{1}{2}\int P \der x} q
\end{align*}
and
\begin{align*}
c_1&=\frac{1}{8}J\chi^2-\frac{1}{8}\left(\frac{1}{\eta}p-\frac{1}{2\eta}(6\xi+P)e^{-\frac{1}{2}\int P \der x} y\right)\chi-\frac{3}{8\eta^3}e^{-\frac{1}{2}\int P \der x}y,\\
c_2&=6z+6\chi J^2-16\left(\frac{1}{\eta}p-\frac{1}{2\eta}(6\xi+P)e^{-\frac{1}{2}\int P \der x} y\right)J\\
&\quad{}-4(6\xi+P)p^2+2(6\xi+P)^2e^{-\frac{1}{2}\int P \der x}yp-\frac{1}{3}(6\xi+P)^3e^{-\int P \der x}y^2\\
&\quad{}+Pp^2+\frac{Q}{10}e^{-\frac{1}{2}\int P \der x}yp+\frac{R}{30}e^{-\int P \der x}y^2,\\
c_3&=-\frac{1}{2}J \chi+\frac{1}{\eta}p-\frac{1}{2\eta}(6\xi+P)e^{-\frac{1}{2}\int P \der x} y,\\
c_4&=-\frac{1}{8}\chi,\\
c_5&=4J.
\end{align*}

Let us write
\begin{align}\label{x1b}
X^1&=\partial_q,\\\label{x2b}
X^2&=\partial_x+e^{\frac{1}{2}\int P \der x}p\partial_y+e^{\frac{1}{2}\int P \der x}q\partial_p+e^{\int P \der x}q^2\partial_z,
\end{align}
and take
\begin{align*}
X^3&=-\eta^3 \chi e^{\frac{1}{2}\int P \der x} \partial_y-\frac{\eta^3 \chi}{2}(6\xi+P)\partial_p+\eta \partial_p\\
&\quad{}+\eta^3 \chi \bigg(e^{-\frac{1}{2}\int P \der x}(12\xi^3+6P\xi^2+(P^2+\frac{Q}{20})\xi+\frac{1}{18}P^3+\frac{PQ}{120}+\frac{R}{90})y\\
&\quad{}-(12\xi^2+3 P \xi+\frac{1}{6}P^2-\frac{Q}{60})p\bigg)\partial_z\\
&\quad{}+\eta(\frac{1}{3}((6\xi+P)^2+\frac{Q}{20})e^{-\frac{1}{2}\int P \der x}y-(8\xi+P)p+4e^{\frac{1}{2}\int P \der x}q)\partial_z.
\end{align*}
This gives 
\begin{align*}
Z^1&=-\eta e^{-\frac{1}{2}\int P dx}\left( \big( 6\xi^2+\frac{3}{2}P\xi+\frac{1}{12}(P^2-\frac{Q}{10})\big)\eta^2 \chi-4\xi-\frac{1}{2} P\right)X^1+X^3,\\
Z^2&=\ell_1X^1-8\eta^2X^2-2\eta\left(\frac{1}{3}((6\xi+P)^2+\frac{Q}{20}) e^{-\frac{1}{2}\int P \der x}y-(8\xi+P) p+2e^{\frac{1}{2}\int P \der x} q\right)X^3,\\
Z^3&=\ell_2X^1-\frac{1}{8}\chi X^3,
\end{align*}
where the functions $\ell_1$ and $\ell_2$ are given by
\begin{align*}
\ell_1&=2e^{-\int P dx} (6\xi^2+\frac{3}{2}P\xi+\frac{1}{12}(P^2-\frac{Q}{10}))\\
&\quad{}\times\left(\frac{1}{3}((6\xi+P)^2+\frac{Q}{20})y-(8\xi+P)e^{\frac{1}{2}\int P \der x}p+2e^{\int P \der x}q\right)\eta^4 \chi \\
&\quad{}-(144\xi^3+68P\xi^2+\frac{32}{3}(P^2+\frac{Q}{32})\xi+\frac{5}{9}(P^3+\frac{9}{100}PQ+\frac{2}{25}R))\eta^2 e^{-\int P \der x}y\\
&\quad{}+(112\xi^2+28 P \xi+\frac{5}{3}P^2-\frac{1}{15}Q)\eta^2 e^{-\frac{1}{2}\int P \der x}p-2(20\xi+P)\eta^2q,\\
\ell_2&=e^{-\frac{1}{2}\int P \der x}\left(\frac{1}{8}(6\xi^2+\frac{3}{2}P\xi+\frac{1}{12}(P^2-\frac{Q}{10}))\eta^3\chi^2-\frac{1}{16}(8\xi+P)\eta \chi+\frac{1}{4\eta}\right).
\end{align*}
We have chosen the preferred ordering of terms according to $\eta^2 \chi$, $e^{-\frac{1}{2}\int P \der x} y$, $p$, $e^{\frac{1}{2}\int P \der x} q$, $\xi$.
We obtain in the end
\begin{align}\label{s1b}
S^1&=-\eta^2\bigg(\frac{2}{9}(216\xi^3+108P\xi^2+18(P^2+\frac{Q}{20})\xi+P^3+\frac{3PQ}{20}+\frac{R}{5})e^{-\int P \der x}y \\
&\quad{}-2(24\xi^2+6P\xi+\frac{1}{3}(P^2-\frac{Q}{10}))e^{-\frac{1}{2}\int P \der x}p+24\xi q\bigg) X^1-8\eta^2X^2,\nonumber \\
\label{s2b}
S^2&=\frac{e^{-\frac{1}{2}\int P \der x}}{4\eta}X^1.
\end{align}
These two vector fields lie in the span of $X^1$ and $X^2$ and together with 
\begin{align}\label{s3b}
S^3&=\ell_3X^1+\bigg((\frac{1}{6}((6\xi+P)^2+\frac{Q}{20})e^{-\frac{1}{2}\int P \der x}y-\frac{1}{2}(8\xi+P)p+e^{\frac{1}{2}\int P \der x}q)\eta^3\chi^2\\
&\quad{}+(\frac{1}{2}(6\xi+P)e^{-\frac{1}{2}\int P \der x}y-p)\eta\chi-\frac{3}{\eta}e^{-\frac{1}{2}\int P \der x}y\bigg)X^2-c_3^2X^3,\nonumber
\end{align} 
where
\begin{align*}
\ell_3&=\eta e^{-\frac{1}{2}\int P dx} \left(\big( 6\xi^2+\frac{3}{2}P\xi+\frac{1}{12}(P^2-\frac{Q}{10})\big)\eta^2 \chi-4\xi-\frac{1}{2}P \right) c_3^2+\frac{e^{-\frac{1}{2}\int  P \der x}}{4\eta}c_2\\
&\quad{}+\eta^2\bigg(\frac{2}{9}(216\xi^3+108P\xi^2+18(P^2+\frac{Q}{20})\xi+P^3+\frac{3PQ}{20}+\frac{R}{5})e^{-\int P \der x}y\\
&\quad{}-2(24\xi^2+6P\xi+\frac{1}{3}(P^2-\frac{Q}{10}))e^{-\frac{1}{2}\int P \der x}p+24\xi q \bigg)c_1,
\end{align*}
they bracket-generate the Lie algebra of split $\frak{g}_2$.

\begin{theorem}\label{g2b}
For the maximally symmetric distribution spanned by $X^1$, $X^{2}$ as given in (\ref{x1b}) and (\ref{x2b}) where $P(x)$ satisfies the generalised Chazy equation with parameter $\frac{2}{3}$ and annihilated by the 1-forms in (\ref{ch231forms}), the Lie algebra of split $\frak{g}_2$ is obtained from the pairwise bracket-generating set $\{S^1,S^2,S^3\}$ where $S^1$, $S^2$ and $S^3$ are given in (\ref{s1b}), (\ref{s2b}), (\ref{s3b}).
\end{theorem}

\section{Example: spin $4$ Lam\'e equation}\label{spin4}
We now give the bracket-generating set of the Lie algebra of split $\frak{g}_2$ when we take solutions of the generalised Chazy equation with $k=\frac{2}{3}$ parametrised by the solutions of the second-order spin $4$ equianharmonic Lam\'e equation
\[
\Phi_{\tilde z \tilde z}-20\wp \Phi=0.
\]
Again take $w=\wp(\tilde z)$ and let $u=4w^3-g_2w-g_3$.
The algebraic form of the Lam\'e equation (with accessory parameter $0$) is given by
\[
u\Phi_{ww}+\frac{1}{2}u_w\Phi_w-n(n+1)w\Phi=0.
\]
When $g_2=0$, $n=4$, we obtain the solutions
\[
\Phi=\al P^{\frac{1}{3}}_{1}\left(\sqrt{1-\frac{4w^3}{g_3}}\right)\sqrt{w}+\beta Q^{\frac{1}{3}}_{1}\left(\sqrt{1-\frac{4w^3}{g_3}}\right)\sqrt{w}
\]
given by associated Legendre functions. They turn out to be algebraic and can be reexpressed as
\[
P^{\frac{1}{3}}_{1}(x)=-\frac{(1-x)^{\frac{5}{6}}}{3(x-1)\Ga(5/3)}(1+x)^{\frac{1}{6}}(3x-1)
\]
and
\[
 Q^{\frac{1}{3}}_{1}(x)=\frac{\pi (1-x)^{\frac{5}{6}}}{3\sqrt{3}(x-1)\Ga(5/3)}\left(\frac{(1-x)^{\frac{1}{3}}(3x+1)}{(1+x)^{\frac{1}{6}}}-\frac{1}{2}(1+x)^{\frac{1}{6}}(3x-1)\right).
\]
Here the branch cut is taken along $(-\infty,-1)\cup(1, \infty)$.
If we introduce the new independent variable $r$ by taking
\[
w=\frac{g_3^{1/3}r}{(r^3+1)^{2/3}},
\]
then the algebraic form of Lam\'e's equation with spin $4$ is 
\[
(r^3+1)^2\Phi_{rr}+2r^2(r^3+1)\Phi_r+20 r \Phi=0,
\]
which has solutions 
\begin{align*}
\Phi=\al\frac{2r^3-1}{(r^3+1)^{4/3}}+\beta \frac{r(r^3-2)}{(r^3+1)^{4/3}}.
\end{align*}
The solution to the spin 4 Lam\'e equation can also be derived from the Hermite-Krichever solution (\cite{krichever} and \cite{maier}), and can also be related to hypergeometric functions and Schwarz triangle functions obtained in \cite{r16} and \cite{r16b}. 

We find that 
\begin{align*}
\der \tilde z=\frac{i}{g_3^{1/6}(r^3+1)^{2/3}}\der r,
\end{align*}
so that taking $k=\frac{2i}{g_3^{1/6}}$ gives
\begin{align*}
\der x=\frac{k}{2(r^3+1)^{2/3}\Phi^2}\der r.
\end{align*}

To determine the solution to (\ref{sode2}), we change the independent coordinate $x$ to $r$ and substitute in $Q=-720\frac{g_3^{1/3}r}{(r^3+1)^{2/3}}\Phi^4$ and find that the general solution is given by
\[
\eta=\frac{\ga (r^3+1)+\delta r (r^3+1)}{\al (2r^3-1)+\beta r (r^3-2)}.
\]
This gives
\begin{align*}
\chi=\int \frac{1}{\eta^2} \der x&=\int \frac{((r^3+1)^{4/3}\Phi)^2}{(r^3+1)^2(\ga+\de r)^2}\frac{k}{2(r^3+1)^{2/3}\Phi^2}\der r\\
&=\int \frac{k}{2(\ga+\de r)^2}\der r\\
&=-\frac{k}{2\de (\de r+\ga)}+c.
\end{align*}
Observe that the spin 1 algebraic form of Lam\'e equation for the same parameter $r$ is
\[
(r^3+1)^2\Phi_{rr}+2r^2(r^3+1)\Phi_r+2 r \Phi=0,
\]
which has solutions 
\begin{align*}
\Phi=\frac{\de r+\ga}{(r^3+1)^{1/3}},
\end{align*}
so $\eta$ can be viewed as a quotient of the spin $4$ solution over the spin $1$ solution, and $\chi$ can be viewed as the independent variable of the generalised Chazy equation with parameter $k=2$ (which corresponds to spin $1$). Again this agrees with the result of Theorem 3.2 in \cite{r19}. 
We shall restrict to the case where $\chi=r$ by taking $\de=0$, $\ga=\sqrt{k/2}$ and the constant of integration $c=0$. This gives us
\begin{align*}
\xi=-\frac{2}{k}\frac{(\beta(r^6+8r^3-2)+9 \al r^2)(\beta(r^4-2r)+\al(2r^3-1))}{(r^3+1)^3}.
\end{align*}
We also have the parametrisation 
\begin{align*}
P&=-\frac{24}{k}\frac{\al(2r^3-1)+\beta r (r^3-2)}{(r^3+1)^3}(\al r^2(r^3-5)-\beta(5r^3-1)),\\
Q&=\frac{2880}{k^2}\frac{r}{(r^3+1)^6}(\al(2r^3-1)+\beta r(r^3-2))^4,\\
R&=\frac{8640}{k^3}\frac{r^3-1}{(r^3+1)^9}(\al(2r^3-1)+\beta r(r^3-2))^6.
\end{align*}

Substituting $F''=\Phi^3$ into (\ref{ch231forms}), the $(2,3,5)$-distribution is in this case given by
\begin{align}\nonumber
\om_1&=\der y-\frac{p}{2}\frac{k \Phi}{(1+r^3)^{2/3}}\der r,\\ \label{lame41form}
 \om_2&=\der p-\frac{q}{2}\frac{k \Phi}{(1+r^3)^{2/3}}\der r,\\ \nonumber
 \om_3&=\der z-\frac{q^2}{2}\frac{k \Phi^4}{(1+r^3)^{2/3}}\der r,
\end{align}
where
\begin{align*}
\Phi=\al\frac{2r^3-1}{(r^3+1)^{4/3}}+\beta \frac{r(r^3-2)}{(r^3+1)^{4/3}} \mbox{~and~} k=\frac{2i}{g_3^{1/6}}.
\end{align*}
Let us take
\begin{align*}
Y=\al(2r^3-1)+\beta r(r^3-2).
\end{align*}

We find
\begin{align}\label{c1b}
c_1&=-\frac{3}{4}\sqrt{\frac{2}{k^3}}y+\frac{\sqrt{2}}{8\sqrt{k}}\frac{r(3\beta r^4+\al (4 r^3+1))}{r^3+1}p+\frac{\sqrt{2k}}{16}\frac{r^2Y^2}{(r^3+1)^3}q,\\ \label{c2b}
c_2&=6z-\frac{288}{k^3}y^2+\frac{8}{k}\frac{(\beta (7 r^4+4r)+\al (8r^3+5))(\beta (2r^3-1)+3\al r^2)}{(r^3+1)^2}p^2\\ \nonumber
&\quad{}+3k\frac{rY^4}{(r^3+1)^6} q^2+\frac{96}{k^2}\frac{(3\beta r^4+\al (4 r^3+1)) r}{r^3+1} yp\\ \nonumber
&\quad{}+\frac{48}{k}\frac{r^2Y^2}{(r^3+1)^3}yq+8\frac{(\al(5r^3+2)+\beta(4r^4+r))Y^2}{(r^3+1)^4}pq,\\ \label{c3b}
c_3&=-\sqrt{\frac{2}{k}}(\beta r+\al)p-\frac{\sqrt{2 k}}{4}\frac{r Y^2}{(r^3+1)^3}q,\\ \nonumber
c_4&=-\frac{r}{8},\\ \nonumber
c_5&=48\sqrt{\frac{2}{k^3}} r y+8\sqrt{\frac{2}{k}}\frac{3\al r^2+\beta(2r^3-1)}{r^3+1}p+2\sqrt{2k}\frac{Y^2}{(r^3+1)^3}q.
\end{align}

This gives
\begin{align*}
\Theta_1&=-\frac{3}{4}\sqrt{\frac{2}{k^3}}\om_1+\frac{1}{8}\sqrt{\frac{2}{k}}\frac{Y}{r^3+1}r\om_2,\\
\Theta_2&=-48\left(\frac{12}{k^3}y+\frac{2(3\al+\beta(r^4+4r))r}{k^2(r^3+1)} p+\frac{r^2Y^2}{k(r^3+1)^3}q\right)\om_1\\
&\quad{}+8\left(\frac{12Yr}{k^2(r^3+1)}y+\frac{2Y(3\al r^2+\beta (2r^3-1))}{k(r^3+1)^2}p-\frac{Y^3}{(r^3+1)^4}q\right)\om_2\\
&\quad{}+6\om_3,\\
\Theta_3&=6\sqrt{\frac{2}{k^3}}r^2\om_1+\sqrt{\frac{2}{k}}\frac{Y}{r^3+1}\om_2.
\end{align*}
Using this set of functions as a basis, we find 
\begin{align*}
Z^1&=\partial_{c_3}+2c_5\partial_{c_2}-2c_4\partial_{c_1}\\
&=-\frac{\sqrt{2}k^{3/2}r}{4(r^3+1)}\partial_y+\frac{2\sqrt{2k}}{(r^3+1)^3}(-\frac{2\beta}{k}(r^3+1)^3p+Y^2q)\partial_z\\
&\quad{}+\sqrt{\frac{k}{2}}\frac{4r^3+1}{Y}\partial_p-2\sqrt{\frac{2}{k}}\frac{(r^3+1)^3(6\al r^2+\beta (5r^3-1))}{Y^3}\partial_q
\end{align*}
and similar formulas for $Z^2$ and $Z^3$.
We have the following corollary:
\begin{corollary}\label{corb}
Consider the maximally symmmetric $(2,3,5)$-distribution given by the span of 
\begin{align*}
X^1=\partial_q, \qquad X^2=\frac{2Y^2}{k(r^3+1)^2}\partial_r+\frac{Y^3}{(r^3+1)^4}p\partial_y+\frac{Y^3}{(r^3+1)^{4}}q\partial_p+\frac{Y^6}{(r^3+1)^{8}}q^2\partial_z
\end{align*}
where $Y=\al(2r^3-1)+\beta r(r^3-2)$.
This is annihilated by the 1-forms $\{\om_1,\om_2,\om_3\}$ in (\ref{lame41form}). Then the vector fields
\begin{align*}
S^1&=24\bigg(\frac{8(r^3+1)^3}{k^2Y^2}y+\frac{4(r^3+1)^2r(\al+\beta r)}{kY^2}p+\frac{9\al r^2+\beta(r^6+8r^3-2)}{(r^3+1)Y}q\bigg)X^1\\
&\quad{}-\frac{4k(r^3+1)^2}{Y^2}X^2,\\
S^2&=\frac{\sqrt{2}(r^3+1)^3}{4\sqrt{k}Y^2}X^1,\\
S^3&=-c_1S^1+c_2S^2-c_3^2Z^1,
\end{align*}
where 
\begin{align*}
Z^1&=-\frac{\sqrt{2}k^{3/2}r}{4(r^3+1)}\partial_y+\frac{2\sqrt{2k}}{(r^3+1)^3}(-\frac{2\beta}{k}(r^3+1)^3p+Y^2q)\partial_z\\
&\quad{}+\sqrt{\frac{k}{2}}\frac{4r^3+1}{Y}\partial_p-2\sqrt{\frac{2}{k}}\frac{(r^3+1)^3(6\al r^2+\beta (5r^3-1))}{Y^3}\partial_q
\end{align*}
and $c_1$, $c_2$ and $c_3$ are given as in (\ref{c1b}), (\ref{c2b}), (\ref{c3b}), pairwise bracket-generate the Lie algebra of split $\frak{g}_2$.
\end{corollary}

The consequence of Theorems \ref{g2a} and \ref{g2b} in this paper is that we obtain additional split $\frak{g}_2$ Lie algbera of vector fields besides those given in (\ref{hcf1}). It would be interesting to investigate the parametrisations of the split real form of the other exceptional Lie algebras ($F_4$, $E_6$, $E_7$ and $E_8$) based on the vector fields that we obtain here. 

In part two of \cite{r18} and later in \cite{r18b}, the automorphisms of the generalised Chazy equation with parameter $k=\frac{3}{2}$ and $k=\frac{2}{3}$ were studied. This gives us additional parametrisations of the $\frak{g}_2$ vector fields based on a particular fixed solution. The transformation to other parameters of the generalised Chazy equation were studied as well, and it would also be interesting to see how it fits into the results obtained in this context.  

\appendix

\section{Split ${\frak g}_2$ Lie algbera of vector fields}
In this appendix we give the split ${\frak g}_2$ vector fields in terms of the coordinate functions $(c_1,c_2,c_3,c_4,c_5)$. 
We have already determined in Section 2 the short roots
\begin{align*}
S^1&=Z^2+c_5Z^1,\\
S^2&=Z^3-c_4Z^1,\\
S^3&=-c_1S^1+c_2S^2-c_3^2Z^1.
\end{align*}
Computing the commutator gives the remaining short roots
\begin{align*}
S^4&=2(\partial_{c_3}+c_4 \partial_{c_1}-c_5\partial_{c_2}),\\
S^5&=2(c_1+3c_3c_4)Z^1-4c_3Z^3+6c_4(c_1\partial_{c_1}+c_2\partial_{c_2}+c_3\partial_{c_3}-\frac{2}{3}c_3\partial_{c_3}-\frac{1}{3}c_4\partial_{c_4}-\frac{1}{3}c_5\partial_{c_5})\\
&\quad{}-6(c_1c_5+c_2c_4+c_3^2)\partial_{c_2},\\
S^6&=-2(c_2-3c_3c_5)Z^1+4c_3Z^2+6c_5(c_1\partial_{c_1}+c_2\partial_{c_2}+c_3\partial_{c_3}-\frac{2}{3}c_3\partial_{c_3}-\frac{1}{3}c_4\partial_{c_4}-\frac{1}{3}c_5\partial_{c_5})\\
&\quad{}-6(c_1c_5+c_2c_4+c_3^2)\partial_{c_1}.
\end{align*}
We obtain the long roots
\begin{align*}
L^1&=-6\partial_{c_1},\\
L^2&=6\partial_{c_2},\\
L^3&=6(c_1\partial_{c_2}-c_4\partial_{c_5}),\\
L^4&=6(c_1-2c_3c_4)((c_1+2c_3c_4)\partial_{c_1}+(c_2-2c_3c_5)\partial_{c_2}+c_3\partial_{c_3}+c_4\partial_{c_4}+c_5\partial_{c_5})\\
&\quad{}+6(c_1c_5+c_2c_4+c_3^2)(2c_3\partial_{c_2}-\partial_{c_5}),\\
L^5&=-6(c_2+2c_3c_5)((c_1+2c_3c_4)\partial_{c_1}+(c_2-2c_3c_5)\partial_{c_2}+c_3\partial_{c_3}+c_4\partial_{c_4}+c_5\partial_{c_5})\\
&\quad{}+6(c_1c_5+c_2c_4+c_3^2)(2c_3\partial_{c_1}+\partial_{c_4}),\\
L^6&=6(c_2\partial_{c_1}-c_5\partial_{c_4}),
\end{align*}
and we also have
\begin{align*}
h&=-6(c_1\partial_{c_1}+c_2\partial_{c_2}+\frac{2}{3}c_3\partial_{c_3}+\frac{1}{3}c_4\partial_{c_4}+\frac{1}{3}c_5\partial_{c_5}),\\
H&=-6(c_2\partial_{c_2}+\frac{1}{3}c_3\partial_{c_3}-\frac{1}{3}c_4\partial_{c_4}+\frac{2}{3}c_5\partial_{c_5}).
\end{align*}
Let $h^1=\frac{1}{4}(h-H)$ and $h^2=\frac{\sqrt{3}}{12}(h+H)$. We have the commutation relations in the table below, where the entries are given by the Lie bracket $[X,Y]$, where $X$ is corresponding vector down the first column and $Y$ is the corresponding vector along the first row. We have displayed only the non-zero elements above the diagonal. (By the anti-symmetry of the Lie bracket, the elements below the diagonal are just the reflection of the upper diagonal entries with the signs reversed).
\begin{table}[h]
\resizebox{\textwidth}{!}{
\bgroup
\def\arraystretch{1.5}
\begin{tabular}{l*{2}{c}|*{6}{c}|*{6}{c}}
$X\backslash Y$      & $h^1$ & $h^2$ & $S^1$ & $S^2$ & $S^3$ & $S^4$ & $S^5$  & $S^6$ & $L^1$ & $L^2$ & $L^3$ & $L^4$ & $L^5$ & $L^6$ \\
\hline
$h^1$                         &    &  & $S^1$ & $-\frac{1}{2}S^2$ & $-\frac{1}{2}S^3$ & $\frac{1}{2}S^4$ & $-S^5$  & $\frac{1}{2}S^6$ & $\frac{3}{2}L^1$ & & $-\frac{3}{2}L^3$ & $\frac{3}{2}L^4$ & & $\frac{3}{2}L^6$ \\
$h^2$                        &     &  &  & $\frac{\sqrt{3}}{2}S^2$ & $-\frac{\sqrt{3}}{2}S^3$ & $\frac{\sqrt{3}}{2}S^4$ &  & $-\frac{\sqrt{3}}{2}S^6$ & $\frac{\sqrt{3}}{2}L^1$ & $\sqrt{3}L^2$ & $\frac{\sqrt{3}}{2}L^3$ & $-\frac{\sqrt{3}}{2}L^4$ & $-\sqrt{3}L^5$ & $-\frac{\sqrt{3}}{2}L^6$ \\
\hline
$S^1$      &  &  &  & $S^4$ & $-S^6$ & $L^1$ & $h-H$  & $L^6$ &  &  & $-6S^2$ & $-6S^3$ &  & \\
$S^2$      &  &  &  &             & $S^5$  & $L^2$ & $L^3$  & $H$ &  &  &  &  & $-6S^3$ & $-6S^1$ \\
$S^3$      &  &  &  &             &              & $-h$ & $L^4$  & $L^5$ & $-6S^1$  & $-6S^2$  &  &  &  &  \\
$S^4$      &  &  &  &             &              &                                        & $-8S^2$  & $8S^1$ &   &  &  & $6S^5$ & $6S^6$ &  \\
$S^5$      &  &  &  &             &              &                                        &   & $-8S^3$ & $6S^4$   &  &  &  & & $6S^6$  \\
$S^6$      &  &  &  &  &  &  & & &  & $6 S^4$ &$6S^5$  &  & &  \\
\hline
$L^1$      &  &  &  &  &  &  & & &  &  &$-6L^2$  &$6(2h-H)$  & $6L^6$ &  \\
$L^2$      &  &  &  &  &  &  & & &  &  &  &$6L^3$  & $6(h+H)$ & $-6L^1$ \\
$L^3$      &  &  &  &  &  &  & & &  &  &  &  & $-6L^4$ & $6(2H-h)$ \\
$L^4$      &  &  &  &  &  &  & & &  &  &  &  &  & $6 L^5$ \\
$L^5$      &  &  &  &  &  &  & & &  &  &  &  &  &  \\
$L^6$      &  &  &  &  &  &  & & &  &  &  &  &  &  
\end{tabular}
\egroup
}
\end{table}

\qr 


\begin{thebibliography}{XX}
\bibitem{acht} M.\ J.\ Ablowitz, S.\ Chakravarty and R.\ Halburd, \textit{The generalized Chazy equation and Schwarzian triangle functions}, Asian J. Math. {\bf 2}, (1998), 1--6.

\bibitem{ach} M.\ J.\ Ablowitz, S.\ Chakravarty and R.\ Halburd, \textit{The generalized Chazy equation from the self-duality equations}, Stud. Appl. Math. {\bf 103}, (1999), 75--88.

\bibitem{AN14} D.\ An and P.\ Nurowski, \textit{Twistor  space  for  rolling  bodies}, 
Comm.  Math.  Phys. {\bf 326}, 2, (2014), 393--414.

\bibitem{annur} D.\ An and P.\ Nurowski, \textit{Symmetric (2,3,5) distributions, an interesting ODE of 7th order and Pleba{\' n}ski metric}, Journ. Geom. Phys. {\bf  126}, (2018), 93--100.

\bibitem{BH} J.\ C.\ Baez and J.\ Huerta, \textit{$G_2$ and the rolling ball}, Trans. Amer. Math. Soc.
{\bf 366} (10), (2014), 5257--5293.

\bibitem{BM} G.\ Bor and R.\ Montgomery, \textit{$G_2$ and the rolling distribution}, Enseign. Math. (2) {\bf 55}, (2009), no. 1-2, 157--196.

\bibitem{cartan1893} E.\ Cartan, \textit{Sur la structure des groupes simples finis et continus}, C.R. Acad. Sc.
{\bf 116} (1893), 784--786.


\bibitem{cartan1910} E.\ Cartan, \textit{Les syst\`emes de Pfaff, \`a cinq variables et les \'equations aux d\'eriv\'ees partielles
du second ordre}, Ann. Sci. \'Ecole Norm. Sup. (3) {\bf 27}
(1910), 109--192.

\bibitem{chazy1} J.\ Chazy, \textit{Sur les \'equations diff\'erentielles dont l'int\'egrale g\'en\'erale est uniforme et admet des singularit\'es essentielles mobiles}, C.R. Acad. Sc. Paris {\bf 149} (1909), 563--565.

\bibitem{chazy2} J.\ Chazy, \textit{Sur les \'equations diff\'erentielles du troisi\`eme ordre et d'ordre sup\'erieur dont l'int\'egrale a ses points critiques fixes}, Acta Math. {\bf 34} (1911), 317--385.

\bibitem{co96} P.\ A.\ Clarkson and P.\ J.\ Olver, \textit{Symmetry and the Chazy Equation}, J. Differential Equations {\bf 124}, 1, (1996), 225--246.

\bibitem{DK} B.\ Doubrov and B.\ Kruglikov, \textit{On the models of submaximal symmetric rank 2 distributions in 5D}, Differential Geom. Appl. {\bf 35} (2014), suppl., 314--322.

\bibitem{Engel} F.\ Engel, \textit{Sur un group simple \`a quatorze param\`etres}, C.R. Acad. Sc. Paris {\bf 116} (1893), 786--788.

\bibitem{Engel2} F.\ Engel, \textit{Zwei merkw\"urdige Gruppen des Raums von f\"unf Dimensionen}, Jahresbericht der Deutschen Mathematiker-Vereinigung {\bf 8} (1900), 196--198.

\bibitem{halphen} G.\ H.\ Halphen, \textit{Sur les invariants des \'equations diff\'erentielles lin\'eaires du quatri\`eme ordre}, Acta Math. {\bf 3} (1883), 325--380. doi:10.1007/BF02422456

\bibitem{krichever} I.\ M.\ Krichever, \textit{Generalized elliptic genera and Baker-Akhiezer functions}, Mat. Zametki {\bf 47} (2), (1990), 34--45, 158. (Transl. Math. Notes {\bf 47}, no. 1-2, (1990), 132--142.)

\bibitem{maier} R.\ S.\ Maier, \textit{Lam\'e polynomials, hyperelliptic reductions and Lam\'e band structure}, Philos. Trans. Roy. Soc. London Ser. A {\bf 366}, (2008), 1115--1153.

\bibitem{conf} P.\ Nurowski, \textit{Differential equations and conformal structures},
Journ. Geom. Phys. {\bf 55}, (2005), 19--49.

\bibitem{r16} M.\ Randall, \textit{Flat (2,3,5)-distributions and Chazy's equations}, SIGMA {\bf 12}, 029, 2016.

\bibitem{r17} M.\ Randall, \textit{SU(2) Pfaffian systems and gauge theory}, arxiv:1705.08172.

\bibitem{r18b} M.\ Randall, \textit{Automorphisms and transformations of solutions to the generalised Chazy equation for various parameters}, J. Differential Equations {\bf 268}, 12, (2020), 7998--8025. 

\bibitem{r16b} M.\ Randall, \textit{Schwarz triangle functions and duality for certain parameters of the generalised Chazy equation}, New Zealand J. Math., {\bf 50}, (2020), 181--205.

\bibitem{r19} M.\ Randall, \textit{Nurowski's conformal class of a maximally symmetric (2,3,5)-distribution and its Ricci-flat representatives}, J. Nonlinear Math. Phys., {\bf 28}, 1, (2021), 1--13.

\bibitem{r18} M.\ Randall, \textit{Automorphism of solutions to Ramanujan's differential equations and other results}, Kyushu J. Math., {\bf 75}, 1, (2021), 77--94.

\bibitem{r21a} M.\ Randall, \textit{A Monge normal form for the rolling distribution}, arxiv:2103.02360.

\bibitem{r21b} M.\ Randall, \textit{Local equivalence of some maximally symmetric $(2,3,5)$-distributions}, arxiv:2108.04599.

\bibitem{Strazzullo} F.\ Strazzullo, \textit{Symmetry Analysis of General Rank-3 Pfaffian Systems in Five Variables},
Ph.D. Thesis, Utah State University (2009).

\bibitem{ww} E.\ T.\ Whittaker and G.\ N.\ Watson, \textit{A Course of Modern Analysis}  (4th ed., Cambridge Mathematical Library). Cambridge: Cambridge University Press, 1996.

\bibitem{tw13} T.\ Willse,
\textit{ Highly symmetric 2-plane fields on 5-manifolds and Heisenberg 5-group holonomy}, Differential Geom. Appl. {\bf 33} (2014), suppl., 81--111.

\bibitem{tw14} T.\ Willse, \textit{Cartan’s incomplete classification and an explicit ambient metric of holonomy $G_2^*$}, Eur. J. Math., {\bf 4}, 2 (2018), 622--638.
\end{thebibliography}
\end{document}